\documentclass[11pt,reqno,oneside]{amsart}

\usepackage{graphicx,subfigure,threeparttable}
\usepackage{amssymb}
\usepackage{epstopdf}
\usepackage{empheq}
\usepackage[a4paper, total={6in, 9in}]{geometry}
\usepackage{graphicx}
\usepackage{float}
\usepackage{subfigure}  
\allowdisplaybreaks[4]

\usepackage{booktabs} 
\usepackage{subfig} 

\usepackage{amsfonts}

\usepackage{amsmath,amsthm,mathrsfs,amssymb,cite,bm}
\usepackage[usenames]{color}

\newtheorem{thm}{Theorem}[section]

\newtheorem{lem}{Lemma}[section]

\theoremstyle{definition}
\newtheorem{defn}{Definition}[section]
\newtheorem{asm}{Assumption}

\newtheorem{rmk}{Remark}

\newcommand{\ii}{\mathrm{i}}
\newcommand{\e}{\mathrm{e}}

\newcommand{\tgao}{\overline{\bm{G}}^{\alpha,\omega}}
\newcommand{\tgaop}[1]{\overline{G}^{\alpha,\omega}_{#1}}
\newcommand{\llm}{\mathcal{L}^{\lambda,\mu}}

\theoremstyle{remark}

\numberwithin{equation}{section}

\numberwithin{equation}{section}

\newcounter{saveeqn}

\newcommand{\doverline}[1]{\overline{\overline{#1}}}

\title[Uniqueness of periodic elastic structures]{Uniqueness for phaseless inverse elastic scattering problem for periodic structures}

\author{Youzi He}
\address{MSU-BIT-SMBU Joint Research Center of Applied Mathematics, Shenzhen MSU–BIT University, 518172 Shenzhen, People's Republic of China}
\address{School of Mathematics and Statistics, Beijing Institute of Technology, 100081 Beijing, People's Republic of China}
\email{youzihe@smbu.edu.cn}
\author{Wei Wu}
\address{School of Mathematics, Jilin University, 130012 Changchun, People's Republic of China}
\email{wei\_wu@jlu.edu.cn}
\thanks{Corresponding author: Wei Wu.}
\author{Hongyi Dang}
\address{School of Mathematics, Jilin University, 130012 Changchun, People's Republic of China}
\email{danghy23@mails.jlu.edu.cn}

\date{} 
\begin{document}

\maketitle

\begin{abstract}
This paper establishes uniqueness results of inverse elastic scattering problem with phaseless near-field data in periodic structures in $\mathbb{R}^2$ and periodic/biperiodic structures in $\mathbb{R}^3$. We use a superposition of two point sources in each periodic unit with different positions as the incident field, and measures the phaseless near-field data on a line parallel to $x_1$-axis in 2D, or on a plane parallel to $(x_1,x_2)$-plane in 3D. We first calculate the explicit formula of quasi-periodic/biperiodic Green's functions of Lam\'{e} system in $\mathbb{R}^3$. Then, to establish the uniqueness results, the reciprocity relations for point sources, scattered fields, and total fields are derived. Finally, with the help of Rayleigh's expansion, the uniqueness results are proved. The quasi-periodic/biperiodic Green's functions of Lam\'{e} system in $\mathbb{R}^3$, the reciprocity relations, and Rayleigh's expansion in $\mathbb{R}^3$ are novel results as important by-products in the proof process.

	\medskip

	\noindent{\bf Keywords:}~~  inverse elastic scattering problem, periodic structures, phaseless inverse scattering problem, uniqueness
	

	
\end{abstract}

\section{Introduction}

The scattering of waves with periodic structures is an important topic in applied mathematics, physics, materials, and engineering. In this broad field, the scattering of elastic waves holds a special position due to its significant industrial applications such as non-destructive testing, seismology, medical imaging, and geophysical prospecting. Unlike acoustic waves or electromagnetic waves, elastic wave propagation is intrinsically more complex, involving the coupling between pressure (P) waves and shear (S) waves. In recent decades, significant progress has been made in the analysis of the forward scattering problem for periodic elastic structures\cite{LiuPRB2000, TrainitiRuzzene}.

However, in many practical scenarios, one is more concerned with the inverse scattering problem. That is, one aims to identify or design the periodic structure from the knowledge of scattered field data. Studying an inverse scattering problem usually needs to address three issues, i.e. existence, uniqueness, and stability. Among these, uniqueness is the theoretical cornerstone of reconstruction algorithms. 

While uniqueness theorems are relatively well-established for inverse scattering problem for periodic structures in the acoustic and electromagnetic cases \cite{ZhangGuo2021phaseless}, corresponding results for the elastic wave system are far less complete. The difficulty in solving this problem comes from the following points. The first is the tensor nature of the elastodynamic system. The second is the intricate coupling of P and S waves. The third is the presence of complex surface waves, which makes it particularly challenging to extract information from the scattering data. Moreover, in most practical applications, the phase information of the waves is difficult to measure. We can only obtain the intensity information of the scattered or total field. This type of inverse problem is known as the phaseless inverse scattering problem. Obviously, compared with the phased case, the phaseless inverse scattering problem is more complicated in the analysis of uniqueness. 

There have been many the existing uniqueness results on phaseless inverse scattering problems for bounded scatterers including \cite{Zhang2010, Zhang2020a, Qu2018, Diao2019, Diao2023, Xu2018, AmmariKang2007PolarizationMomentTensors}, and research papers such as \cite{NiuLvGao, Ammari1995UniquenessDoublyPeriodic} demonstrated uniqueness results on periodic structures in acoustic and electromagnetic scenario. In this work, we consider the inverse elastic scattering problem for periodic structures with phaseless data. We consider the scattering problem of elastic waves in 2D-periodic structure and 3D-periodic and biperiodic structures. In 2-dimensional case, the structure is described by a $C^2$ function $x_2=f_1(x_1)$ with periodicity 1 in $x_1$. In 3-dimensional cases, two structures are described by $C^2$ functions $x_3=f_2(x_1)$ and $x_3=f_3(x_1,x_2)$, where $f_2$ has periodicity 1, and $f_3$ is periodic with respect to both $x_1$ and $x_2$ with periodicity $1$ in both directions. We denote those two structures as periodic and biperiodic structures, correspondingly. We use a superposition of two point sources in each periodic unit with different positions as the incident field. To explicitly illustrate the scattered wave, we first calculate the quasi-periodic and quasi-biperiodic Green's functions of Lam\'{e} system in $\mathbb{R}^3$, which are also important by-products as they are the basis of constructing solutions to arbitrary source problems. We directly cite the results from \cite{WuHe} for quasi-periodic Green's function in $\mathbb{R}^2$. For periodic structures, since the scattered field contains the evanescent waves, it is impossible to establish the one-to-one correspondence between the scattered fields and the far field patterns. Therefore, we turn to deriving the reciprocity relations for point sources, scattered fields, and total fields, respectively. These relations, together with Rayleigh's expansion, are necessary to obtain the desired uniqueness results. The derivation of Rayleigh's expansion of quasi-periodic system in $\mathbb{R}^3$ is another important contribution of this work.

The rest of this paper is organized as follows. In section~\ref{sec2}, we introduce the mathematical model of the scattering problem for both the quasi-periodic system in $\mathbb{R}^2$ and the quasi-periodic/biperiodic system in $\mathbb{R}^3$. We also deduce the Rayleigh expansion of scattered wave of quasi-periodic system in $\mathbb{R}^3$ in this section. In section~\ref{sec3}, we calculate the quasi-periodic and quasi-biperiodic Green's functions of Lam\'{e} system in $\mathbb{R}^3$. The reciprocity relations for point sources, scattered fields, and total fields are given in section~\ref{sec4}. Section~\ref{sec5} is committed to achieving the uniqueness results of the phaseless inverse elastic scattering problems for periodic structures. The auxiliary theorem needed in the proof of final theorems is concluded in the Appendix with their proofs.

\section{Formulation of Problems}\label{sec2}

\subsection{Formulation of problem in $\mathbb{R}^2$}

In this section, we consider the scattering of elastic waves by an infinite periodic structure in $\mathbb{R}^2$. Assume that the smooth surface $\Gamma_1$ is described by a $C^2$-smooth function with periodicity 1, i.e.,
\begin{equation}\label{eq:Gamma}
\Gamma_1 =\{(x_1,x_2)\in \mathbb{R}^2: x_2=f(x_1), f(x_1+1) = f(x_1) \}.
\end{equation}
The surface $\Gamma_1$ is considered to be impenetrable and displacement free, i.e. the displacement field is 0 on $\Gamma_1$. Denote the unit normal and tangential vectors by $\bm{\nu}$ and $\bm{\tau}$, respectively. Denote by $\mathbb{S}:=\{\bm{x}\in \mathbb{R}^2 : |\bm{x}|=1\}$ the unit circle. The whole space is divided by $\Gamma_1$ into two parts. In a single periodic unit, we denote the region above $\Gamma_1$ by
\[
\Omega_1^{+} = \{(x_1,x_2)\in\mathbb{R}^2: x_2>f(x_1), 0\leq x_1\leq 1\}.
\]
We collect measured data on the line
\[
\Gamma_{h,1} =\{(x_1,h): x_1\in\mathbb{R} \},
\]
where $h$ is a fixed constant satisfying $h> \max_{x_1\in[0,1]}f(x_1)$. Denote the domain between $\Gamma_1$ and $\Gamma_{h,1}$ in one period by
\[
D_1 = \{(x_1,x_2)\in\mathbb{R}^2: 0\leq x_1\leq 1, f(x_1)<x_2<h \}.
\]
Denote the left boundary of $D_1$ as $\Gamma_{L,1}:=\{(x_1,x_2)\in\mathbb{R}^2:x_1=0, f(x_1)<x_2<h \}$, the right boundary of $D_1$ as $\Gamma_{R,1}:=\{(x_1,x_2)\in\mathbb{R}^2:x_1=1, f(x_1)<x_2<h \}$. We assume the region $D_1$ is occupied by an isotropic, homogeneous elastic medium characterized by the Lam\'e constants $\lambda, \mu\ (\mu>0, \lambda+ \mu>0)$ and density $\rho=1$. In such system, the wavenumber of p-wave is $k_p:=\omega\sqrt{1/(\lambda+2\mu)}$, and the wavenumber of s-wave is $k_s=\omega\sqrt{1/\mu}$, where $\omega>0$ is the angular frequency. The unit $D_1$ is periodically aligned along $x_1$ direction with periodicity 1.

\begin{figure}
\centering
\includegraphics[width=0.6\textwidth]{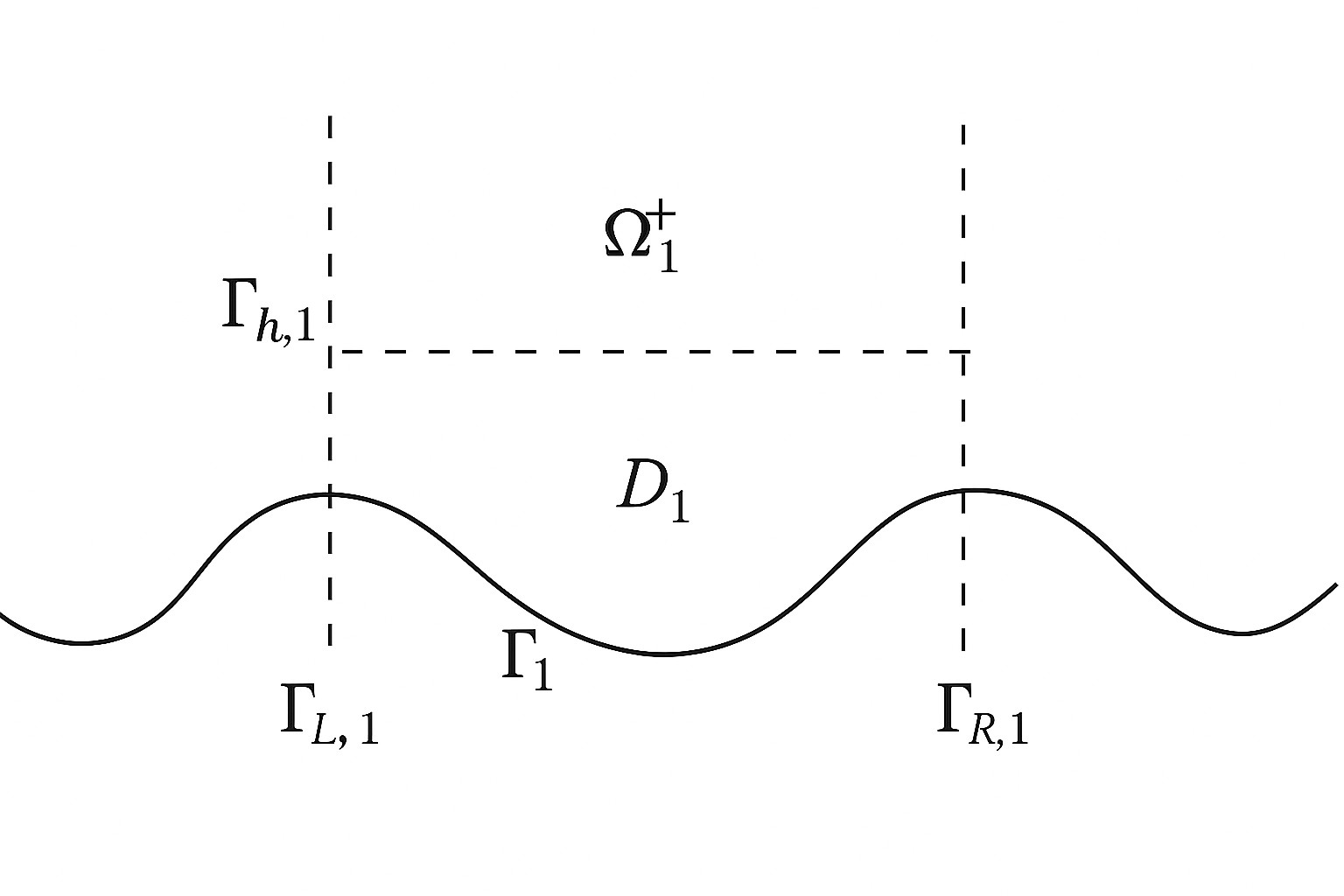}
\caption{Quasi-periodic system in $\mathbb{R}^2$}
\end{figure}
The diffraction of the incident plane wave $\bm{u}_{{\rm{inc}}}(\bm{x})$ illuminating $\Gamma_1$ can be formulated to the time-reduced Navier equation
\begin{equation}\label{eq:Navier}
\Delta^{\ast}\bm{u}_{{\rm{tot}}}+ \omega^2\bm{u}_{{\rm{tot}}} =\bm{0} \quad \text{in }\ \Omega_1^{+},
\end{equation}
with the boundary condition
\begin{equation}\label{eq:boundaryC}
\bm{u}_{{\rm{tot}}} =\bm{0} \quad \text{on }\ \Gamma_1.
\end{equation}
Here, $\Delta^{\ast}:=\mu \Delta + (\lambda +\mu)\nabla\nabla\cdot$, $\bm{u}_{{\rm{tot}}}=\bm{u}_{{\rm{inc}}} + \bm{u}_{{\rm{sc}}}$ is the total field, $\bm{u}_{{\rm{sc}}}$ is the scattered field. Due to the periodicity of the whole system, it is enough to calculate in a single periodic unit.



The incident field is composed of waves from different point sources. To illustrate the position of point sources, we introduce the following definition of admissible curve as \cite{ZGSL, ChengDong}.
\begin{defn}\label{def:admissiblecurve}
An open set $\Sigma_1$ is called an \emph{admissible curve} with respect to simply-connected domain $\Omega_1$ if
\begin{enumerate}
\item  $\Omega_1\subset D_1$ bounded and simply connected;
\item  $\partial \Omega_1$ is analytic homeomorphic to the unit circle;
\item  $\omega^2$ is not a Dirichlet eigenvalue of $-\Delta^{\ast}$ in $\Omega_1$;
\item  $\Sigma_1\subset \partial \Omega_1$ is a one-dimensional analytic manifold with non-vanishing measure.
\end{enumerate}
\end{defn}

In the current system, we use a superposition of two point sources in each periodic unit with different positions as the incident field. Let $\Gamma_{0,1}$ be a line parallel to $\Gamma_{h,1}$ and positioned between $\Gamma_1$ and $\Gamma_{h,1}$. Choose a fixed point $\tilde{\bm{z}}$ on $\Gamma_{0,1}$, where a point source with polarization vector $\bm{q}$ is placed. Next, in each periodic unit we select an admissible curve $\Sigma_1 \subset\partial \Omega_1$, where $\Omega_1$ is an open set between $\Gamma_{h,1}$ and $\Gamma_{0,1}$. We place another point source with polarization vector $\bm{q}'$ on $\bm{z}\in \Sigma_1$.

Denote $\bm{G}(\bm{x},\bm{y})$ be the free-space Green's function of Lam\'{e} system satisfying
$$
    (\Delta^*_x + \omega^2)\bm{G}(\bm{x},\bm{y}) = \delta_{\bm{y}}(\bm{x}).
$$

The radiating elastic wave from a single point source at $\bm{y}$ with polarization $\bm{q}$ is $\bm{G}(\bm{x},\bm{y})\bm{q}$. We require the neighboring point sources in the system having a phase difference of $\alpha$. Let $\bm{G}^\alpha(\bm{x},\bm{y})$ denote the quasi-periodic Green's function of Lam\'{e} system, which will be discussed in detail in Section \ref{sec:green2d}. The superposition of all point sources located at $\tilde{\bm{z}}+(n,0)$, $n\in\mathbb{Z}$ is $\bm{G}^\alpha(\bm{x}, \widetilde{\bm{z}})\bm{q}$, and the superposition of all point sources located at $\bm{z}+(n,0)$, $n\in\mathbb{Z}$ is $\bm{G}^\alpha(\bm{x},\bm{z})\bm{q}'$.

Let $\bm{u}_{\alpha,{\rm{sc}}}(\bm{x},\bm{z},\bm{q}')$ (resp. $\bm{u}_{\alpha,{\rm{tot}}}(\bm{x},\bm{z},\bm{q}')$) and $\bm{u}_{\alpha,{\rm{sc}}}(\bm{x},\tilde{\bm{z}},\bm{q})$ (resp. $\bm{u}_{\alpha,{\rm{tot}}}(\bm{x},\tilde{\bm{z}},\bm{q})$) be the scattered fields (resp. total fields) generated by $\bm{G}^{\alpha}(\bm{x},\bm{z})\bm{q}'$ and $\bm{G}^{\alpha}(\bm{x},\tilde{\bm{z}})\bm{q}$, respectively.

Obviously, the incident field and the scattered field contain outgoing waves in the region above the point sources. To avoid the cancellation of the outgoing waves of the incident field and the scattered field, we need the following assumption.

\begin{asm}\label{asm:1}
For any periodic structure, the total field generated by the point sources contains at least one upgoing wave.
\end{asm}

It is shown in \cite{Arens99} that the scattered field $\bm{u}_{{\rm{sc}}}(\bm{x})$ is also a quasi-periodic solution with momentum $\alpha$. According to \cite{Arens99, Charalambopoulos01}, for $x_2>\max_{x_1\in[0,1]}f(x_1)$, the scattered field admits the following expansion
\begin{align}\label{eq:sc_radiation}
\bm{u}_{{\rm{sc}}}(\bm{x}) = \sum_{l\in 2\pi\mathbb{Z}}&u_l^p(\alpha_l\hat{\bm{e}}_1+\beta_l\hat{\bm{e}}_2 ) {\rm{e}}^{{\rm{i}} (\alpha_l x_1+\beta_l x_2)} +\sum_{l\in 2\pi\mathbb{Z}}u_l^s (\gamma_l\hat{\bm{e}}_1-\alpha_l \hat{\bm{e}}_2) {\rm{e}}^{{\rm{i}} (\alpha_l x_1+\gamma_l x_2)},
\end{align}
where $u_l^p$ and $u_l^s$ are arbitrary complex constants and
\begin{align*}
\alpha_l := \alpha + l, && \beta_l := \left\{
\begin{aligned}
    \sqrt{k_p^2-\alpha_l^2}, && \alpha_l^2\leq k_p^2,\\
    \mathrm{i}\sqrt{\alpha_l^2-k_p^2}, && \alpha_l^2>k_p^2,
\end{aligned}\right. &&
 \gamma_l := \left\{
 \begin{aligned}
    \sqrt{k_s^2-\alpha_l^2}, && \alpha_l^2\leq k_s^2,\\
    \mathrm{i}\sqrt{\alpha_l^2-k_s^2}, && \alpha_l^2>k_s^2.
 \end{aligned}\right.
\end{align*}

\subsection{Formulation of problem in $\mathbb{R}^3$}

Now we introduce the quasi-periodic/biperiodic system in $\mathbb{R}^3$.

\begin{figure}
\centering
\includegraphics[width=0.6\textwidth]{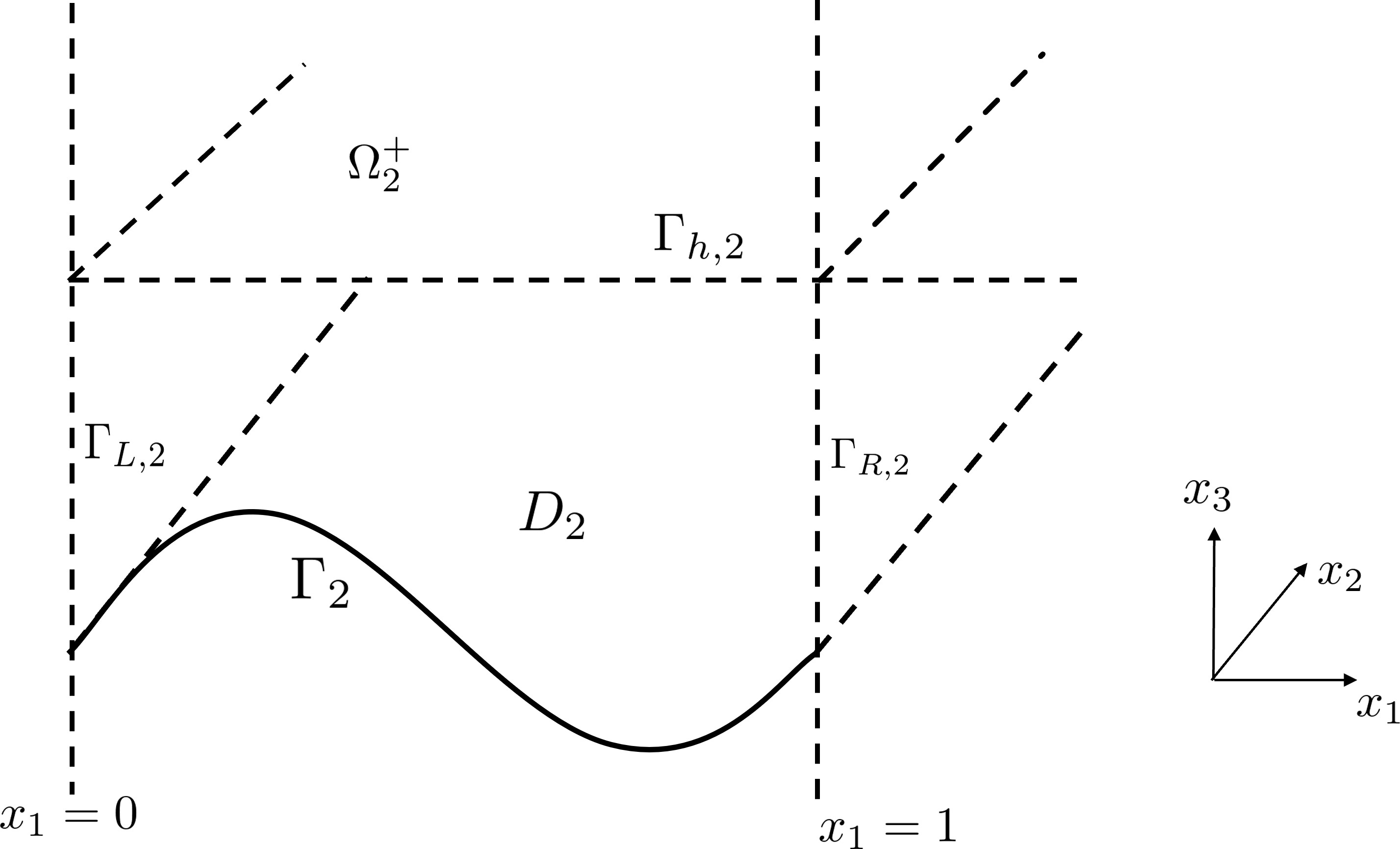}
\caption{Quasi-periodic system in $\mathbb{R}^3$}
\label{fig:3d1}
\end{figure}

\begin{figure}
\centering
\includegraphics[width=0.6\textwidth]{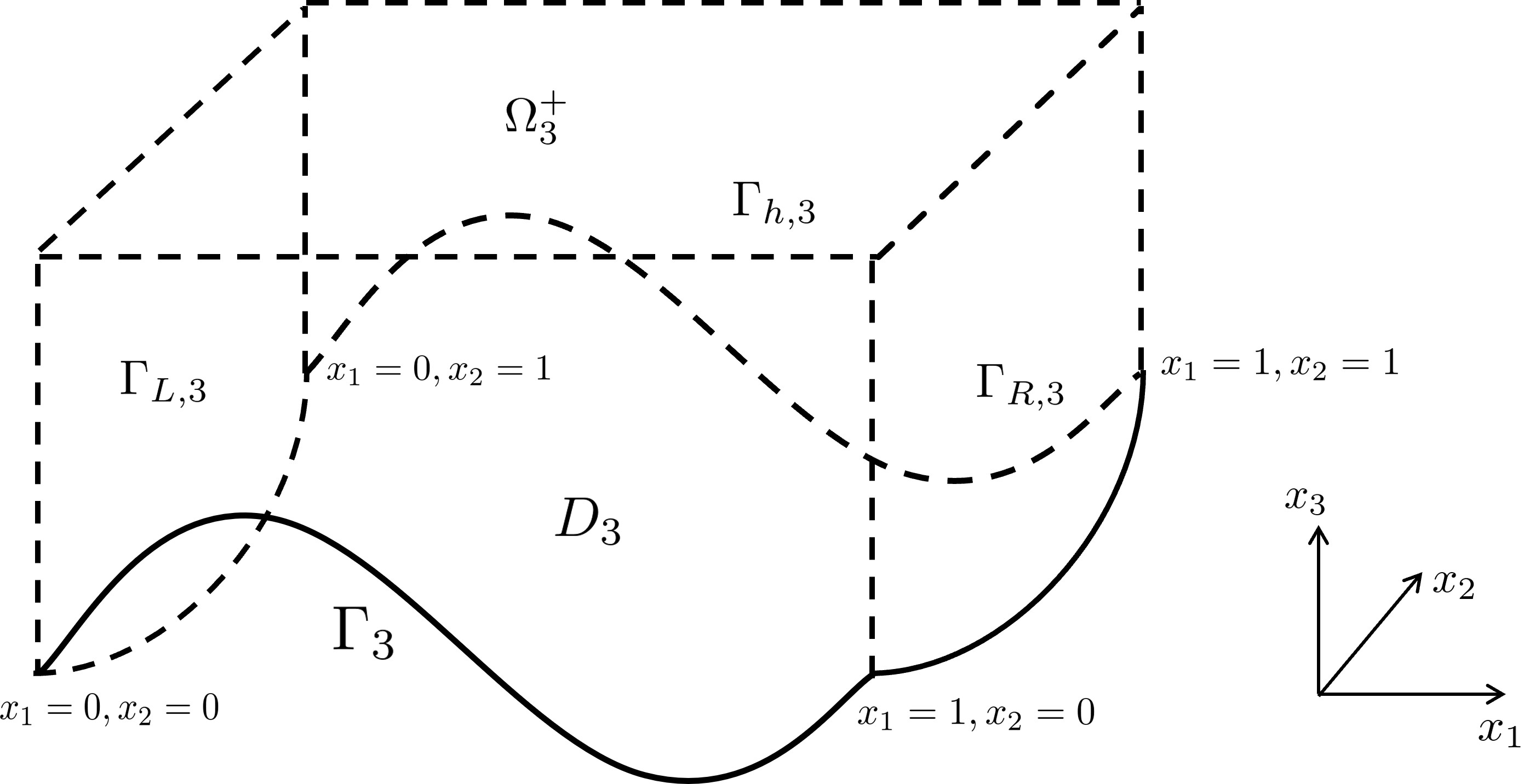}
\caption{Quasi-biperiodic system in $\mathbb{R}^3$}
\label{fig:3d2}
\end{figure}

Assume that the smooth surface $\Gamma_2, \Gamma_3$ are described as follows
\begin{align*}
\Gamma_2 &=\{(x_1,x_2,x_3)\in\mathbb{R}^3: x_3=f_2(x_1), f_2(x_1+1) = f_2(x_1) \}, \\
\Gamma_3 &=\{(x_1,x_2,x_3)\in\mathbb{R}^3: x_3=f_3(x_1,x_2), f_3(x_1+1,x_2) = f_3(x_1,x_2), f_3(x_1, x_2+1) = f_3(x_1,x_2)\}.
\end{align*}
The surface $\Gamma_2, \Gamma_3$ are also considered to be impenetrable and displacement free. Denote by $\mathbb{S}^2:=\{\bm{x}\in \mathbb{R}^3 : |\bm{x}|=1\}$ the unit sphere. $\Omega_i^+, \Gamma_{h,i}, D_i, \Gamma_{L,i}, \Gamma_{R,i}, i=2,3$ demonstrated in Figure \ref{fig:3d1} and Figure \ref{fig:3d2} are defined in the similar way
\begin{align*}
    \Omega_2^+ &:= \{\bm{x}\in\mathbb{R}^3: x_3>f_2(x_1), 0\leq x_1\leq 1, x_2\in\mathbb{R}\}, \\
    \Omega_3^+ &:= \{\bm{x}\in\mathbb{R}^3: x_3>f_3(x_1,x_2), 0\leq x_1, x_2 \leq 1\}, \\
    \Gamma_{h,i} &:= \{(x_1, x_2, h_i): x_1, x_2\in\mathbb{R} \},~~ h_2>\max_{x_1\in[0,1]}f_i(x_1), ~~ h_3>\max_{(x_1,x_2)\in[0,1]^2}f(x_1,x_2), \\
    D_2 &:= \{\bm{x}\in\mathbb{R}^3: 0 \leq x_1 \leq 1, f_2(x_1)<x_3<h_2, x_2\in\mathbb{R}\}, \\
    D_3 &:= \{\bm{x}\in\mathbb{R}^3: 0\leq x_1, x_2\leq 1, f_3(x_1,x_2)<x_3<h_3 \}, \\
    \Gamma_{L,2}&:= \{\bm{x}\in\mathbb{R}^3: x_1=0, f_2(x_1)<x_3<h_2, x_2\in\mathbb{R}\}, \\
    \Gamma_{R,2}&:= \{\bm{x}\in\mathbb{R}^3: x_1=1, f_2(x_1)<x_3<h_2, x_2\in\mathbb{R}\}, \\
    \Gamma_{L,3}&:= \{\bm{x}\in\mathbb{R}^3: x_1=0, 0 \leq x_2 \leq 1, f_3(x_1,x_2)<x_3<h_3 \}, \\
    \Gamma_{R,3}&:= \{\bm{x}\in\mathbb{R}^3: x_1=1, 0 \leq x_2 \leq 1, f_3(x_1,x_2)<x_3<h_3 \}.
\end{align*}
We assume both quasi-periodic and quasi-biperiodic systems are occupied by an isotropic, homogeneous elastic medium with Lam\'e constants $\lambda, \mu\ (\mu>0, \lambda+ \mu>0)$ and density $\rho=1$. In quasi-periodic system, the unit $D_2$ is periodically aligned along $x_1$ direction with periodicity 1, while in quasi-biperiodic system, the unit $D_3$ is periodically aligned along $x_1$ and $x_2$ direction with periodicity 1.

\vspace{1em}

Different from 2-dimensional case, here we will instead use the following definition of admissible surface as \cite{ZhangGuo2021phaseless}.
\begin{defn}
An open set $\Sigma_i$ is called an \emph{admissible surface} with respect to simply-connected domain $\Omega_i$ if
\begin{enumerate}
\item  $\Omega_i\subset D_i$ bounded and simply connected;
\item  $\partial \Omega_i$ is analytic homeomorphic to the unit sphere;
\item  $\omega^2$ is not a Dirichlet eigenvalue of $-\Delta^{\ast}$ in $\Omega_i$;
\item  $\Sigma_i\subset \partial \Omega_i$ is a two-dimensional analytic manifold with non-vanishing measure.
\end{enumerate}
\end{defn}

We still use a superposition of two point sources in each periodic unit with different positions as the incident field. Let $\Gamma_{0,i}, i=2,3$ be a plane parallel to $\Gamma_{h,i}$ and positioned between $\Gamma_i$ and $\Gamma_{h,i}$. Choose a fixed point $\tilde{\bm{z}}_i$ on $\Gamma_{0,i}$, where a point source with polarization vector $\bm{q}$ is placed. Next, in each periodic unit we select an admissible surface $\Sigma_i \subset\partial \Omega_i$, where $\Omega_i$ is an open set between $\Gamma_{h,i}$ and $\Gamma_{0,i}$. We place another point source with polarization vector $\bm{q}'$ at $\bm{z}_i\in \Sigma_i$.

In quasi-periodic system, we denote $\alpha\in\mathbb{R}_+$ as quasi-periodicity, which means point sources with position differed by $(1,0,0)$ have phase difference $\alpha$. In quasi-biperiodic system, we denote $\bm{\alpha}:=(\alpha_1,\alpha_2)$ as the quasi-periodicity. Sources with position difference $(1,0,0)$ have a phase difference $\alpha_1$, and sources with position difference $(0,1,0)$ have a phase difference $\alpha_2$. Let $\widetilde{\bm{G}}^\alpha(\bm{x},\bm{y})$ and $\mathring{\bm{G}}^{\bm{\alpha}}(\bm{x},\bm{y})$ denote the quasi-periodic and quasi-biperiodic Green's function in Lam\'{e} system in $\mathbb{R}^3$ respectively, which will be calculated in detail in Section \ref{sec:3d1} and \ref{sec:3d2}. The superposition of all point sources located at $\tilde{\bm{z}}+(n,0,0)$, $n\in\mathbb{Z}$(resp. $\tilde{\bm{z}}+(n,m,0)$, $n,m\in\mathbb{Z}$ is $\widetilde{\bm{G}}^\alpha(\bm{x}, \widetilde{\bm{z}})\bm{q}$(resp. $\mathring{\bm{G}}^{\bm{\alpha}}(\bm{x}, \widetilde{\bm{z}})\bm{q}$), and the superposition of all point sources located at $\bm{z}+(n,0,0)$, $n\in\mathbb{Z}$(resp. $\bm{z}+(n,m,0)$, $n,m\in\mathbb{Z}$) is $\mathring{\bm{G}}^{\bm{\alpha}}(\bm{x},\bm{z})\bm{q}'$.

When there is no ambiguity, we still use $\bm{u}_{\alpha,{\rm{sc}}}(\bm{x},\bm{z},\bm{q}')$ (resp. $\bm{u}_{\alpha,{\rm{tot}}}(\bm{x},\bm{z},\bm{q}')$) and $\bm{u}_{\alpha,{\rm{sc}}}(\bm{x},\tilde{\bm{z}},\bm{q})$ (resp. $\bm{u}_{\alpha,{\rm{tot}}}(\bm{x},\tilde{\bm{z}},\bm{q})$) be the scattered fields (resp. total fields) generated by $\widetilde{\bm{G}}^{\alpha}(\bm{x},\bm{z})\bm{q}'$ (or $\mathring{\bm{G}}^{\bm{\alpha}}(\bm{x},\bm{z})\bm{q}'$) and $\widetilde{\bm{G}}^{\alpha}(\bm{x},\tilde{\bm{z}})\bm{q}$(or $\mathring{\bm{G}}^{\bm{\alpha}}(\bm{x},\tilde{\bm{z}})\bm{q}$), respectively. To ensure the existence of outgoing wave, we still insist Assumption \ref{asm:1} as before.

\subsection{Rayleigh expansion of quasi-periodic/biperiodic system in $\mathbb{R}^3$}

We firstly deduce the Rayleigh expansion of scattered wave of quasi-periodic system in $\mathbb{R}^3$. In $\Omega_2^+$, the scattered wave in elastic system could be decomposed as (see e.g. \cite{HuKirschSini2013})
$$
    \bm{u}^{\mathrm{sc}} = \nabla\phi + \nabla\times\bm{\psi}
$$
with $\phi = -\frac{1}{k_p^2}\nabla\cdot \bm{u}^{\mathrm{sc}}$, $\bm{\psi} = \frac{1}{k_s^2}\nabla\times \bm{u}^{\mathrm{sc}}$.
In $\Omega$, the scalar function $\phi$ satisfies $(\Delta + k_p^2)\phi = 0$, while the vector function $\bm{\psi}$ satisfies $(\Delta + k_s^2)\bm{\psi} = 0$. Both $\phi$ and $\bm{\psi}$ satisfies the quasi-periodic condition \eqref{eq:quasi_period_2}, which means we have the following expansion with respect to $x_1$
\begin{align*}
&\phi(\bm{x}) = \sum\limits_{\substack{n\in 2\pi\mathbb{Z} \\ (\alpha+n)^2\leq k_p^2}} \phi_n(x_2, x_3) \e^{\ii(\alpha+n)x_1}, \quad \bm{\psi}(\bm{x}) = \sum\limits_{\substack{n\in 2\pi\mathbb{Z} \\ (\alpha+n)^2\leq k_s^2}} \bm{\psi}_n(x_2,x_3) \e^{\ii(\alpha+n)x_1}.
\end{align*}
The additional requirements imposed on $n$ is because that only outgoing solutions are considered in the scattered wave, therefore we need $(\alpha+n)^2<k_p^2$ to ensure propagating p-wave exists, and $(\alpha+n)^2<k_s^2$ to ensure propagating s-wave exists. By plugging the expansion into the Helmholtz equation they satisfy, we could immediately conclude that for $\phi_n, \bm{\psi}_n$,
\begin{align*}
    \Delta \phi_n + (k_p^2-(\alpha+n)^2)\phi_n = 0, \quad \Delta \bm{\psi}_n + (k_s^2-(\alpha+n)^2)\bm{\psi}_n = 0.
\end{align*}
Then, $\phi_n, \bm{\psi}_n$ takes the form of
$$
    \phi_n(x_2,x_3) = \sum\limits_{m\in\mathbb{Z}} A_m H_m^{(1)}(\beta_n r)\e^{\ii m\theta}, \quad \bm{\psi}_n(x_2,x_3) = \sum\limits_{m\in\mathbb{Z}} \bm{B}_m H_m^{(1)}(\gamma_n r)\e^{\ii m\theta},
$$
where $r:=\sqrt{x_2^2+x_3^2}, \tan\theta=x_3/x_2, \beta_n = \sqrt{k_p^2-(\alpha+n)^2}, \gamma_n = \sqrt{k_s^2-(\alpha+n)^2}$, and $A_m\in\mathbb{C}, \bm{B}_m:=(B_{m,1}, B_{m,2}, B_{m,3})\in\mathbb{C}^3$ be constants. Define $\widetilde{\bm{\psi}}_n(x_1,x_2,x_3):=\bm{\psi}_n(x_2,x_3)$, and
\begin{align*}
    \bm{d}_{p,n} := \nabla\phi_n = \sum\limits_{m\in\mathbb{Z}} A_m\frac{\e^{\ii m\theta}}{r}&\left(x_2\beta_n(H_m^{(1)}(\beta_n r))' - \frac{\ii mx_3}{r}H_m^{(1)}(\beta_n r),\right. \\& \left.x_3\beta_n (H_m^{(1)}(\beta_n r))' + \frac{\ii mx_2}{r}H_m^{(1)}(\beta_n r)\right)
\end{align*}
\begin{equation*}
    \bm{d}_{s,n} := \nabla\times\widetilde{\bm{\psi}}_n := \left((\nabla\times\widetilde{\bm{\psi}}_n)_1,(\nabla\times\widetilde{\bm{\psi}}_n)_2,(\nabla\times\widetilde{\bm{\psi}}_n)_3\right),
\end{equation*}
where
\begin{align*}
&(\nabla\times\widetilde{\bm{\psi}}_n)_1 = \\&\sum\limits_{m\in\mathbb{Z}}\frac{\e^{\ii m\theta}}{r}\left[\left(x_2\gamma_n(H_m^{(1)}(\gamma_n r))'-\frac{\ii mx_3}{r}H_m^{(1)}(\gamma_n r)\right)B_{m,3} - \left(x_3\gamma_n(H_m^{(1)}(\gamma_n r))'+\frac{\ii mx_2}{r}H_m^{(1)}(\gamma_n r)\right)B_{m,2}\right],\\
&(\nabla\times\widetilde{\bm{\psi}}_n)_2 = \sum\limits_{m\in\mathbb{Z}}\frac{\e^{\ii m\theta}}{r}B_{m,1}\left(x_3\gamma_n(H_m^{(1)}(\gamma_n r))' + \frac{\ii mx_2}{r}H_m^{(1)}(\gamma_n r)\right),\\
&(\nabla\times\widetilde{\bm{\psi}}_n)_3 = \sum\limits_{m\in\mathbb{Z}}\frac{\e^{\ii m\theta}}{r}B_{m,1}\left(\frac{\ii m x_3}{r}H_m^{(1)}(\gamma_n r) - x_2\gamma_n(H_m^{(1)}(\gamma_n r))'\right).
\end{align*}
The scattered field then admits the following expansion (Rayleigh expansion)
\begin{equation}\label{eq:rayleigh31}
\begin{aligned}
\bm{u}^{\mathrm{sc}}(\bm{x}) &= \sum\limits_{\substack{n\in 2\pi\mathbb{Z} \\ (\alpha+n)^2\leq k_p^2}}\e^{\ii(\alpha+n)x_1}\left(\ii(\alpha+n)\phi_n, \bm{d}_{p,n}\right) \\ &+  \sum\limits_{\substack{n\in 2\pi\mathbb{Z} \\ (\alpha+n)^2\leq k_s^2}}\e^{\ii(\alpha+n)x_1}\left(\bm{d}_{s,n} + \ii(\alpha+n)\sum\limits_{m\in\mathbb{Z}}\e^{\ii m\theta}H_m^{(1)}(\gamma_n r)(0, -B_{m,3}, B_{m,2})\right).
\end{aligned}
\end{equation}

Here we directly quote results from \cite{ElschnerHu2012} on the Rayleigh expansion of quasi-biperiodic elastic system in $\mathbb{R}^3$. Suppose the quasi-biperiodicity is $\bm{\alpha}:=(\alpha_1,\alpha_2)$, and define $\bar{\bm{x}}:=(x_1,x_2)$. Then for $x_3$ large enough, it holds that
\begin{equation}\label{eq:rayleigh32}
\begin{aligned}
u^{\mathrm{sc}}(\bm{x}) = \sum\limits_{\bm{n}\in\mathbb{Z}^2} A_{p,n}(\bm{\alpha}_n,\beta_n)^\top \exp(\ii\bm{\alpha}_n\cdot\bar{x} + \ii\beta_nx_3) +\bm{A}_{s,n}\exp(\ii\bm{\alpha}_n\cdot\bar{x} + \ii\gamma_nx_3).
\end{aligned}
\end{equation}
Here $A_{p,n}\in\mathbb{C}, \bm{A}_{s,n}\in\mathbb{C}^3$ are constants, $\bm{\alpha}_n:=(\alpha_1+2\pi n_1, \alpha_2+2\pi n_2)$, $n_1,n_2\in\mathbb{Z}$, and $\beta_n, \gamma_n$ are defined by
\begin{equation*}
\beta_n = \left\{
\begin{aligned}
\sqrt{k_p^2-|\bm{\alpha}_n|^2} && \mathrm{if}~|\bm{\alpha}_n|\leq k_p, \\
\ii\sqrt{|\bm{\alpha}_n|^2-k_p^2} && \mathrm{if}~|\bm{\alpha}_n|>k_p,
\end{aligned}
\right.
\quad\gamma_n = \left\{\begin{aligned}
\sqrt{k_s^2-|\bm{\alpha}_n|^2} && \mathrm{if}~|\bm{\alpha}_n|\leq k_s, \\
\ii\sqrt{|\bm{\alpha}_n|^2-k_s^2} && \mathrm{if}~|\bm{\alpha}_n|>k_s.
\end{aligned}\right.
\end{equation*}
It is worth mentioning that \eqref{eq:rayleigh31} does not include evanescent wave in $x_3$ direction, while \eqref{eq:rayleigh32} includes the evanescent wave in $x_3$ direction. Since our discussion only concerns the propagating wave in $x_3$ direction, the difference between \eqref{eq:rayleigh31} and \eqref{eq:rayleigh32} is negligible.

\section{Quasi-periodic Green's functions in 2D and 3D}\label{sec3}

In this section we demonstrate the calculation of quasi-periodic and quasi-biperiodic Green's functions of Lam\'{e} equation in 3D. For the sake of completeness, we will also show the expression of quasi-periodic Green's function in 2D, which was already calculated in \cite{WuHe}.

\begin{defn}
Function $\bm{u}_1:\mathbb{R}^2\rightarrow\mathbb{R}^2$ is called $\alpha$-quasi-periodic (in $x_1$) with momentum $\alpha>0$ if the vector function $\bm{v}_1(\bm{x}) = {\rm{e}}^{-{\rm{i}}\alpha x_1}\bm{u}_1(\bm{x})$ is periodic with respect to $x_1$, for all $\bm{r}\in\mathbb{R}^2$, that is the quasi-periodic function $\bm{u}$ satisfies
\begin{equation}\label{eq:quasi_period}
\bm{u}_1(x_1+n, x_2) = {\rm{e}}^{{\rm{i}}n\alpha}\bm{u}_1(x_1, x_2), n\in\mathbb{Z}.
\end{equation}
Similarly, function $\bm{u}_2:\mathbb{R}^3\rightarrow\mathbb{R}^3$ is called $\alpha$-quasi-periodic (in $x_1$) with momentum $\alpha>0$ if
\begin{equation}\label{eq:quasi_period_2}
    \bm{u}_2(x_1+n, x_2, x_3) = \e^{\ii n\alpha}\bm{u}_2(x_1,x_2,x_3), n\in\mathbb{Z}.
\end{equation}
Function $\bm{u}_3: \mathbb{R}^3\rightarrow\mathbb{R}^3$ is called $\bm{\alpha}$-quasi-biperiodic (in $x_1$ and $x_2$) with momentum $\bm{\alpha}:=(\alpha_1, \alpha_2)$, $\alpha_1, \alpha_2>0$ if
\begin{equation}\label{eq:quasi_period_3}
    \bm{u}_3(x_1+n_1, x_2+n_2, x_3) = \e^{\ii(n_1\alpha_1+n_2\alpha_2)}\bm{u}_3(x_1,x_2,x_3), n_1, n_2\in\mathbb{Z}.
\end{equation}
The quasi-periodicity is actually the part of wave vector projected onto the direction of periodicity. Hence $|\alpha|$ in quasi-periodic system or $|\bm{\alpha}|$ in quasi-biperiodic system can never be larger than p-wave number $k_p$ and s-wave number $k_s$. 
\end{defn}
\begin{rmk}
In this paper, we assume all quasi-periodic Green's functions in 2D and 3D are quasi-periodic in $x_1$, and all quasi-biperiodic Green's functions in 3D are quasi-biperiodic in $x_1$ and $x_2$. For brevity we will not mention the periodic variable in the remaining of this paper.
\end{rmk}

\subsection{2D quasi-periodic Green's function}\label{sec:green2d}
The two-dimensional $\alpha$-quasi-periodic Green's function $\bm{G}^{\alpha,\omega}$ for the Lam\'e system satisfies

\begin{equation}\label{eq:quasi-p-Green}
\begin{split}
(\mathcal{L}^{\lambda, \mu} +\rho\omega^2\bm{I})\bm{G}^{\alpha,\omega}(\bm{x},\bm{y})&=\sum_{n \in \mathbb{Z}} \delta(x_1 - y_1 - n) \delta(x_2 - y_2) {\rm{e}}^{{\rm{i}} n\alpha}\bm{I}\\
&=\delta(x_2 - y_2)\left(\sum_{n \in \mathbb{Z}} \delta(x_1 - y_1 - n){\rm{e}}^{{\rm{i}} n\alpha} \right)\bm{I}.
\end{split}
\end{equation}

From \cite{WuHe}, we know that it can be represented by
\begin{small}
\begin{align*}
& \bm{G}^{\alpha,\omega}(\bm{x},\bm{y})= \sum_{l\in L_1}\frac{{\rm{i}}}{4\pi}{\rm{e}}^{\mathrm{i}\alpha_l(x_1-y_1)} \frac{\lambda+\mu}{\mu(\lambda+2\mu)(k_p^2-k_s^2)}\\
  &\bigg(
  \begin{array}{ll}
   \sqrt{k_s^2- \alpha_l^2}{\rm{e}}^{\mathrm{i}\sqrt{k_s^2- \alpha_l^2}|x_2-y_2|}+\frac{\alpha_l^2}{\sqrt{k_p^2-\alpha_l^2}} {\rm{e}}^{\mathrm{i}\sqrt{k_p^2- \alpha_l^2}|x_2-y_2|}    &
      \text{sgn}(x_2-y_2)\alpha_l ({\rm{e}}^{\mathrm{i}\sqrt{k_p^2- \alpha_l^2}|x_2-y_2|}-{\rm{e}}^{\mathrm{i}\sqrt{k_s^2- \alpha_l^2}|x_2-y_2|} )   \\
    \text{sgn}(x_2-y_2)\alpha_l ({\rm{e}}^{\mathrm{i}\sqrt{k_p^2- \alpha_l^2}|x_2-y_2|}-{\rm{e}}^{\mathrm{i}\sqrt{k_s^2- \alpha_l^2}|x_2-y_2|} )    &
      \sqrt{k_p^2- \alpha_l^2}{\rm{e}}^{\mathrm{i}\sqrt{k_p^2- \alpha_l^2}|x_2-y_2|}+\frac{\alpha_l^2}{\sqrt{k_s^2- \alpha_l^2}}{\rm{e}}^{\mathrm{i}\sqrt{k_s^2- \alpha_l^2}|x_2-y_2|}
       \end{array}
       \bigg)\\
       &+\sum_{l\in L_2}\frac{1}{4\pi}{\rm{e}}^{\mathrm{i}\alpha_l(x_1-y_1)} \frac{\lambda+\mu}{\mu(\lambda+2\mu)(k_s^2-k_p^2)}\\
&\bigg(
  \begin{array}{ll}
   -\frac{\alpha_l^2}{\sqrt{ \alpha_l^2- k_p^2}}{\rm{e}}^{-\sqrt{\alpha_l^2-k_p^2}|x_2-y_2|}-\mathrm{i} \sqrt{k_s^2- \alpha_l^2}{\rm{e}}^{\mathrm{i}\sqrt{k_s^2- \alpha_l^2}|x_2-y_2|}  &
     \mathrm{i}\alpha_l \text{sgn}(x_2-y_2) ({\rm{e}}^{\mathrm{i}\sqrt{k_s^2- \alpha_l^2}|x_2-y_2|}-{\rm{e}}^{-\sqrt{ \alpha_l^2-k_p^2}|x_2-y_2|} )   \\
    \mathrm{i}\alpha_l \text{sgn}(x_2-y_2) ({\rm{e}}^{\mathrm{i}\sqrt{k_s^2- \alpha_l^2}|x_2-y_2|}-{\rm{e}}^{-\sqrt{ \alpha_l^2-k_p^2}|x_2-y_2|} )   &
      \sqrt{\alpha_l^2-k_p^2}{\rm{e}}^{-\sqrt{\alpha_l^2-k_p^2}|x_2-y_2|} - \frac{\mathrm{i}\alpha_l^2}{\sqrt{k_s^2- \alpha_l^2}}{\rm{e}}^{\mathrm{i}\sqrt{k_s^2- \alpha_l^2}|x_2-y_2|}
       \end{array}
       \bigg)\\
       &+\sum_{l\in L_3}\frac{1}{4\pi}{\rm{e}}^{\mathrm{i}\alpha_l(x_1-y_1)} \frac{\lambda+\mu}{\mu(\lambda+2\mu)(k_s^2-k_p^2)}\\
&\bigg(
  \begin{array}{ll}
  \sqrt{\alpha_l^2-k_s^2}{\rm{e}}^{-\sqrt{\alpha_l^2-k_s^2}|x_2-y_2|}-\frac{\alpha_l^2}{\sqrt{\alpha_l^2-k_p^2}} {\rm{e}}^{-\sqrt{\alpha_l^2-k_p^2}|x_2-y_2|}    &
   \mathrm{i}\alpha_l \text{sgn}(x_2-y_2) ({\rm{e}}^{-\sqrt{\alpha_l^2-k_s^2}|x_2-y_2|}-{\rm{e}}^{-\sqrt{\alpha_l^2-k_p^2}|x_2-y_2|} )   \\
   \mathrm{i}\alpha_l \text{sgn}(x_2-y_2) ({\rm{e}}^{-\sqrt{\alpha_l^2-k_s^2}|x_2-y_2|}-{\rm{e}}^{-\sqrt{\alpha_l^2-k_p^2}|x_2-y_2|} )  &
   \sqrt{\alpha_l^2-k_p^2}{\rm{e}}^{-\sqrt{\alpha_l^2-k_p^2}|x_2-y_2|} - \frac{\alpha_l^2}{\sqrt{\alpha_l^2-k_s^2}}{\rm{e}}^{-\sqrt{\alpha_l^2-k_s^2}|x_2-y_2|}
    \end{array}
       \bigg) \\
   &=: \sum\limits_{l\in L_1}{\rm{e}}^{\mathrm{i}\alpha_l(x_1-y_1)}\bm{G}^{\alpha_l}_1(x_2,y_2) + \sum\limits_{l\in L_2}{\rm{e}}^{\mathrm{i}\alpha_l(x_1-y_1)} \bm{G}^{\alpha_l}_2(x_2,y_2) + \sum\limits_{l\in L_3}{\rm{e}}^{\mathrm{i}\alpha_l(x_1-y_1)}\bm{G}^{\alpha_l}_3(x_2,y_2).
\end{align*}
\end{small}
where $\alpha_l:=\alpha+l$, and
\begin{equation}\label{eq:defl123}
\begin{aligned}
&L_1:=\{l\in 2\pi\mathbb{Z}:\alpha_l^2<k_p^2 \},\\
&L_2:=\{l\in 2\pi\mathbb{Z}:k_p^2\leq \alpha_l^2<k_s^2 \}, \\
&L_3:=\{l\in 2\pi\mathbb{Z}:\alpha_l^2\geq k_s^2 \}.
\end{aligned}
\end{equation}
\begin{rmk}\label{rem:2d}
Notice that, in the expression of $\bm{G}^{\alpha_l}_i$, $\alpha_l$ is either paired with $\mathrm{sgn}(x_2-y_2)$, or in the form of $\alpha_l^2$ together with $|x_2-y_2|$. Therefore, $\bm{G}^{-\alpha_l}_i(y_2,x_2) = \bm{G}^{\alpha_l}_i(x_2,y_2)$, $i=1,2,3$.
\end{rmk}

\subsection{3D quasi-periodic Green's function}\label{sec:3d1}
The three-dimensional $\alpha$-quasi-periodic Green's function of Lam\'{e} system  $\widetilde{\bm{G}}^{\alpha,\omega}$ satisfies
\begin{equation}\label{eq:quasi-p-Green-3d}
\begin{split}
(\mathcal{L}^{\lambda, \mu} +\rho\omega^2\bm{I}_3)\widetilde{\bm{G}}^{\alpha,\omega}(\bm{x},\bm{y})&=\sum_{n \in \mathbb{Z}} \delta(x_1 - y_1 - n) \delta(x_2 - y_2)\delta(x_3-y_3) {\rm{e}}^{{\rm{i}} n\alpha}\bm{I}_3\\
&=\delta(x_2 - y_2)\delta(x_3-y_3)\left(\sum_{n \in \mathbb{Z}} \delta(x_1 - y_1 - n){\rm{e}}^{{\rm{i}} n\alpha} \right)\bm{I}_3.
\end{split}
\end{equation}
Here $\bm{I}_3$ is $3\times 3$ identity matrix, and $\widetilde{\bm{G}}^{\alpha,\omega}(\bm{x},\bm{y})$ is $3\times 3$ matrix-valued function. It is straightforward that $\widetilde{\bm{G}}^{\alpha,\omega}(\bm{x},\bm{y}) = \widetilde{\bm{G}}^{\alpha,\omega}(\bm{x}-\bm{y})$. Define $\overline{\bm{G}}^{\alpha,\omega}(\bm{x}):= e^{-\ii \alpha x_1}\widetilde{\bm{G}}^{\alpha,\omega}(\bm{x})$. $\tgao(\bm{x})$ is then a periodic function of $x_1$ with periodicity 1. Substituting $e^{\ii \alpha x_1}\overline{\bm{G}}$ into \eqref{eq:quasi-p-Green-3d} we get the following results.
    \begin{equation}
        \begin{aligned}
            e^{-\ii\alpha x_1}(\llm(e^{\ii\alpha x_1} \tgao))_{1j} &= \mu(-\alpha^2+2\ii\alpha\partial_1 + \Delta)\tgaop{1j}+(\lambda+\mu)(-\alpha^2+2\ii\alpha\partial_1+\partial_{11})\tgaop{1j} \\ &+ (\lambda+\mu)(\ii\alpha\partial_2+\partial_{12})\tgaop{2j} + (\lambda+\mu)(\ii\alpha\partial_3 + \partial_{13})\tgaop{3j}, \\
            e^{-\ii\alpha x_1}(\llm(e^{\ii\alpha x_1}\tgao))_{kj} &= \mu(-\alpha^2+2\ii\alpha\partial_1+\Delta)\tgaop{kj}+(\lambda+\mu)(\ii\alpha\partial_k+\partial_{1k})\tgaop{1j} \\ &+(\lambda+\mu)(\partial_{2k}\tgaop{2j}+\partial_{3k}\tgaop{3j}).
        \end{aligned}
    \end{equation}
Here $j=1,2,3, k=2,3$, and $\partial_i f$ is short for $\frac{\partial f}{\partial x_i}$, $\partial_{il} f$ is short for $\frac{\partial^2 f}{\partial x_i \partial x_l}$. Define $\overline{\mathcal{L}}^{\lambda,\mu} := (\overline{\mathcal{L}}^{\lambda,\mu}_{ij})_{3\times 3}, i,j=1,2,3$ as
    \begin{align*}
    \overline{\mathcal{L}}^{\lambda,\mu}_{11} &= (\lambda+2\mu)(-\alpha^2+2\ii\alpha\partial_1+\partial_{11}) + \mu(\partial_{22} + \partial_{33}), \\
    \overline{\mathcal{L}}^{\lambda,\mu}_{ii} &= \mu(-\alpha^2+2\ii\alpha\partial_1+\Delta) + (\lambda+\mu)\partial_{ii}, \quad i=2,3,\\
    \overline{\mathcal{L}}^{\lambda,\mu}_{12} &= \overline{\mathcal{L}}^{\lambda,\mu}_{21} = (\lambda+\mu)(\ii\alpha\partial_2 + \partial_{12}), \\
    \overline{\mathcal{L}}^{\lambda,\mu}_{13} &= \overline{\mathcal{L}}^{\lambda,\mu}_{31} = (\lambda+\mu)(\ii\alpha\partial_3 + \partial_{13}), \\
    \overline{\mathcal{L}}^{\lambda,\mu}_{23} &= \overline{\mathcal{L}}^{\lambda,\mu}_{32} = (\lambda+\mu)\partial_{23}.
    \end{align*}
Then
\begin{equation}\label{eq:tgao1}
    (\overline{\mathcal{L}}^{\lambda,\mu}+\rho\omega^2\bm{I})\tgao(\bm{x}) = \sum\limits_{n\in\mathbb{Z}}\delta(x_1-n)\delta(x_2)\delta(x_3)\bm{I}_3.
\end{equation}
Consider the Fourier expansion of $\tgao(\bm{x})$ with respect to $x_1$ in the form of
\begin{equation}\label{eq:tgaopij}
    \tgaop{ij}(x_1,x_2,x_3) = \sum\limits_{l\in 2\pi\mathbb{Z}} c_{l,ij}(x_2,x_3)e^{\ii lx_1}.
\end{equation}
To get rid of $x_1$ dependency in the equation, we take advantage of Poisson summation formula
$$
    \sum\limits_{n\in\mathbb{Z}}\delta(x_1-n) = \frac{1}{2\pi}\sum\limits_{l\in 2\pi\mathbb{Z}}e^{\ii lx_1}
$$
Plugging Fourier expansion and Poisson summation formula into \eqref{eq:tgao1} leads to an ordinary differential equation system of $\bm{c}_l(x_2,x_3):=(c_{l,ij}(x_2,x_3))_{3\times 3}$.The expressions are shown below.

\begin{lem}\label{eq:cl}
    $\bm{c}_l$ satisfies
    \begin{equation}\label{eq:clode}
        \begin{aligned}
            &(-(\lambda+2\mu)(\alpha+l)^2+\rho\omega^2)c_{l,1j} + \mu (\partial_{22}+\partial_{33})c_{l,1j} + \ii(\lambda+\mu)(\alpha+l)(\partial_2 c_{l,2j} + \partial_3 c_{l,3j}) = \frac{1}{2\pi}\delta_0\delta_{1j}, \\
            &\ii(\lambda+\mu)(\alpha+l)\partial_2c_{l,1j} + (-\mu(\alpha+l)^2+\rho\omega^2)c_{l,2j} + (\lambda+2\mu)\partial_{22}c_{l,2j} + \mu\partial_{33}c_{l,3j} + (\lambda+\mu)\partial_{23}c_{l,3j} = \frac{1}{2\pi}\delta_0\delta_{2j}, \\
            &\ii(\lambda+\mu)(\alpha+l)\partial_3c_{l,1j} + (\lambda+\mu)\partial_{23}c_{l,2j} + (-\mu(\alpha+l)^2+\rho\omega^2)c_{l,3j} + \mu\partial_{22}c_{l,3j} + (\lambda+2\mu)\partial_{33}c_{l,3j} = \frac{1}{2\pi}\delta_0\delta_{3j}.
        \end{aligned}
    \end{equation}
    for $j=1,2,3$.
\end{lem}

Notice that $c_{l,1j}, c_{l,2j}, c_{l,3j}$ are paired in the equation system for a fixed $j=1,2,3$. Applying Fourier transform to \eqref{eq:clode} and letting $|\bm{\xi}|^2:=\xi_2^2+\xi_3^2$, $r:=\sqrt{x_2^2+x_3^2}$, we get
\begin{equation*}\label{eq:clft}
\left\{
    \begin{aligned}
        &(\rho\omega^2-(\lambda+2\mu)(\alpha+l)^2-\mu|\xi|^2)\hat{c}_{l,1j} - (\lambda+\mu)(\alpha+l)\xi_2\hat{c}_{l,2j}-(\lambda+\mu)(\alpha+l)\xi_3\hat{c}_{l,3j}= (2\pi)^{-2}\delta_{1j}, \\
        &-(\lambda+\mu)(\alpha+l)\xi_2\hat{c}_{l,1j} + (-\mu(\alpha+l)^2+\rho\omega^2-(\lambda+2\mu)\xi_2^2)\hat{c}_{l,2j} - \mu\xi_3^2\hat{c}_{l,2j} - (\lambda+\mu)\xi_2\xi_3\hat{c}_{l,3j} = (2\pi)^{-2}\delta_{2j}, \\
        &-(\lambda+\mu)(\alpha+l)\xi_3\hat{c}_{l,1j} - (\lambda+\mu)\xi_2\xi_3\hat{c}_{l,2j} + (-\mu(\alpha+l)^2+\rho\omega^2-(\lambda+2\mu)\xi_3^2)\hat{c}_{l,3j} - \mu\xi_2^2\hat{c}_{l,3j} =(2\pi)^{-2}\delta_{3j}.
    \end{aligned}
\right.
\end{equation*}
Define
\begin{align*}
D &= (\rho\omega^2-(\lambda+2\mu)(\alpha+l)^2-\mu|\bm{\xi}|^2)(\rho\omega^2-\mu(\alpha+l)^2-(\lambda+2\mu)|\bm{\xi}|^2)-(\lambda+\mu)^2(\alpha+l)^2|\bm{\xi}|^2 \\
&= \mu(\lambda+2\mu)(|\bm{\xi}|^2+(\alpha+l)^2-\frac{\rho\omega^2}{\lambda+2\mu})(|\bm{\xi}|^2+(\alpha+l)^2-\frac{\rho\omega^2}{\mu}).
\end{align*}

The straightforward calculation shows that
\begin{equation}\label{eq:hatc}
\begin{aligned}
\hat{c}_{l,11}(\bm{\xi}) &= \frac{1}{4\pi^2}  \frac{\rho\omega^2-\mu(\alpha+l)^2-(\lambda+2\mu)|\bm{\xi}|^2}{D},\\
\hat{c}_{l,21}(\bm{\xi})
&= \hat{c}_{l,12}(\bm{\xi}) = \frac{1}{4\pi^2}\frac{(\lambda+\mu)(\alpha+l)\xi_2}{D}, \\
\hat{c}_{l,31}(\bm{\xi})
&= \hat{c}_{l,13}(\bm{\xi}) = \frac{1}{4\pi^2}\frac{(\lambda+\mu)(\alpha+l)\xi_3}{D}, \\
\hat{c}_{l,22}(\bm{\xi}) &= \frac{1}{4\pi^2}\frac{\rho\omega^2-(\lambda+2\mu)(\alpha+l)^2 -\mu\xi_2^2 - (\lambda+2\mu)\xi_3^2}{D}, \\
\hat{c}_{l,32}(\bm{\xi}) & = \hat{c}_{l,23}(\bm{\xi}) = \frac{1}{4\pi^2}\frac{(\lambda+\mu)\xi_2\xi_3}{D}, \\
\hat{c}_{l,33}(\bm{\xi}) &= \frac{1}{4\pi^2}\frac{\rho\omega^2-(\lambda+2\mu)(\alpha+l)^2 - (\lambda+2\mu)\xi_2^2-\mu\xi_3^2 }{D}.
\end{aligned}
\end{equation}





To calculate the inverse Fourier transforms, we quote necessary results in the following lemma.

\begin{lem}[]\label{lem:equality1}
From (6.566.2) of \cite{gradshteyn2014}, the following integration holds
    \begin{equation}\label{eq:jmukmu}
        \int_0^{+\infty} \frac{x^{\mu+1}J_\mu(ax)}{b^2+x^2}\mathrm{d}x = b^\mu K_\mu(ab)
    \end{equation}
    for $a>0,Re(b)>0, -1<Re(\mu)<3/2$.

From (9.6.26) of \cite{abramowitz1964handbook}, modified Bessel function $K_\nu(z)$ satisfies the following recursive equation
\begin{equation}\label{eq:recursivek}
    K_\nu^\prime(z) = -K_{\nu-1}(z) - \frac{\nu}{z}K_\nu(z).
\end{equation}

From (9.1.41) of \cite{abramowitz1964handbook}, we have
\begin{equation}\label{eq:generatingfunction}
    e^{\frac{1}{2}z(t-\frac{1}{t})} = \sum\limits_{n=-\infty}^{+\infty}J_n(z)t^n.
\end{equation}
\end{lem}

From Lemma \ref{lem:equality1} we get the following useful results.
\begin{lem}\label{lem:aux1}
\begin{align*}
    &\frac{1}{2\pi}\int_{\mathbb{R}^2}\frac{\e^{\ii\bm{\xi}\cdot\bm{x}}}{h^2+|\bm{\xi}|^2}\mathrm{d}\bm{\xi} = K_0(hr), \\
    &\frac{1}{2\pi}\int_{\mathbb{R}^2}\frac{\xi_2}{h^2+|\bm{\xi}|^2}\e^{\ii\bm{\xi}\cdot\bm{x}}\mathrm{d}\bm{\xi} = \frac{\ii hx_2}{r} K_1(hr), \\
    &\frac{1}{2\pi}\int_{\mathbb{R}^2}\frac{\xi_2\xi_3}{h^2+|\bm{\xi}|^2}\e^{\ii\bm{\xi}\cdot\bm{x}}\mathrm{d}\bm{\xi} = -\frac{x_2x_3h}{r^2}\left(\frac{2}{r}K_1(hr) + hK_0(hr)\right), \\
    &\frac{1}{2\pi}\int_{\mathbb{R}^2}\frac{\xi_2^2}{h^2+|\bm{\xi}|^2}\e^{\ii\bm{\xi}\cdot\bm{x}}\mathrm{d}\bm{\xi} = \frac{h}{r^2}\left(\frac{x_3^2-x_2^2}{r}K_1(hr) - x_2^2 hK_0(hr)\right),
\end{align*}
where $r=|\bm{x}|$.
\end{lem}
\begin{proof}
  Here we only prove the first equality. The other equalities follow by taking derivatives of the first equality.

  If we let $t=\e^{\ii\phi}$ in \eqref{eq:generatingfunction},
  $$
    e^{\frac{1}{2}z(t-\frac{1}{t})} = e^{\ii z\sin\phi} = e^{\ii z\cos(\phi-\frac{\pi}{2})} = \sum\limits_{n=-\infty}^{+\infty} J_n(z)\e^{\ii n\phi}.
  $$
  Therefore
  \begin{equation}\label{eq:ge1}
    \e^{\ii z\cos\phi} = \sum\limits_{n=-\infty}^{+\infty} J_n(z) \e^{\ii n(\phi+\frac{\pi}{2})} = \sum\limits_{n=-\infty}^{+\infty} \ii^n J_n(z)e^{\ii n\phi}.
  \end{equation}
  Let $\xi_2 = r_\xi\cos\theta$, $\xi_3=r_\xi\sin\theta$, $x_2=r\cos\phi$, $x_3=r\sin\phi$, then
  \begin{align*}
    \int_{\mathbb{R}^2}\frac{\e^{\ii\bm{\xi}\cdot\bm{x}}}{h^2+|\bm{\xi}|^2}\mathrm{d}\bm{\xi} &= \int_0^{2\pi}\int_0^{+\infty}\frac{\e^{\ii r_\xi r\cos(\theta-\phi)}}{h^2+r_\xi^2}r_\xi\mathrm{d}r_\xi\mathrm{d}\theta = \int_0^{+\infty}\frac{r_\xi}{h^2+r_\xi^2}\int_0^{2\pi}\e^{\ii r_\xi r\cos(\theta-\phi)}\mathrm{d}\theta\mathrm{d}r_\xi \\
    &= \int_0^{+\infty}\frac{r_\xi}{h^2+r_\xi^2}\int_0^{2\pi}\sum\limits_{n=-\infty}^{+\infty}\ii^n J_n(r_\xi r)\e^{\ii n(\theta-\phi)}\mathrm{d}\theta\mathrm{d}r_\xi = \int_0^{+\infty} \frac{2\pi r_\xi J_0(r_\xi r)}{h^2+r_\xi^2}\mathrm{d}r_\xi \\
    &= 2\pi K_0(hr).
  \end{align*}
  The last equality holds by using \eqref{eq:jmukmu} with $\mu=0, a=r, b=h$.
\end{proof}

\begin{lem}\label{lem:aux2}
\begin{align*}
    &\frac{1}{2\pi}\int_{\mathbb{R}^2}\frac{\e^{\ii\bm{\xi}\cdot\bm{x}}}{-h^2+|\bm{\xi}|^2}\mathrm{d}\bm{\xi} = -\frac{\pi\ii}{2}H_0^{(1)}(hr), \\
    &\frac{1}{2\pi}\int_{\mathbb{R}^2}\frac{\xi_2}{-h^2+|\bm{\xi}|^2}\e^{\ii\bm{\xi}\cdot\bm{x}}\mathrm{d}\bm{\xi} = \frac{\pi hx_2}{2r} H_1^{(1)}(hr), \\
    &\frac{1}{2\pi}\int_{\mathbb{R}^2}\frac{\xi_2\xi_3}{-h^2+|\bm{\xi}|^2}\e^{\ii\bm{\xi}\cdot\bm{x}}\mathrm{d}\bm{\xi} = -\frac{\pi\ii x_2x_3h}{2r}\left(-\frac{2}{r}H_1^{(1)}(hr) + hH_0^{(1)}(hr)\right), \\
    &\frac{1}{2\pi}\int_{\mathbb{R}^2}\frac{\xi_2^2}{-h^2+|\bm{\xi}|^2}\e^{\ii\bm{\xi}\cdot\bm{x}}\mathrm{d}\bm{\xi} = \frac{\pi\ii h}{2r^2}\left(\frac{x_3^2-x_2^2}{r}H_1^{(1)}(hr) - x_2^2 hH_0^{(1)}(hr)\right),
\end{align*}
where $r=|\bm{x}|$.
\end{lem}
\begin{proof}
The proof is similar to Lemma \ref{lem:aux1}.
\end{proof}

Now we calculate the explicit expressions of $c_{l,ij}$. Define $\gamma_{l,s} := \sqrt{|-(\alpha+l)^2+k_s^2|}, \gamma_{l,p} := \sqrt{|-(\alpha+l)^2+k_p^2|}$. We divide the calculations into three different cases depending on $l$:
\begin{itemize}
\item \textbf{Case I: $(\alpha+l)^2\geq k_s^2>k_p^2$.}

In this case $D = \mu(\lambda+2\mu)(|\bm{\xi}|^2+\gamma_{l,s}^2)(|\bm{\xi}|^2+\gamma_{l,p}^2)$. Using Lemma \ref{lem:aux1} and with a bit complicated calculations, we have
\begin{align*}
    c_{l,11}(x_2,x_3) &= \frac{1}{4\pi^2\rho\omega^2}\gamma_{l,s}^2K_0(\gamma_{l,s}r) - \frac{(\alpha+l)^2}{4\pi^2\rho\omega^2}K_0(\gamma_{l,p}r), \\
    c_{l,21}(x_2,x_3) &= c_{l,12}(x_2,x_3) = \frac{\ii(\alpha+l)x_2}{4\pi^2\rho\omega^2 r}(\gamma_{l,s}K_1(\gamma_{l,s}r) - \gamma_{l,p}K_1(\gamma_{l,p}r)), \\
    c_{l,31}(x_2,x_3) &= c_{l,13}(x_2,x_3) = \frac{\ii(\alpha+l)x_3}{4\pi^2\rho\omega^2 r}(\gamma_{l,s}K_1(\gamma_{l,s}r) - \gamma_{l,p}K_1(\gamma_{l,p}r)) , \\
    c_{l,22}(x_2,x_3) &= -\frac{1}{4\pi^2\mu}K_0(\gamma_{l,s}r) + \frac{1}{4\pi^2\rho\omega^2r^2}\left(\gamma_{l,s}\left(\frac{x_3^2-x_2^2}{r}K_1(\gamma_{l,s}r) - x_2^2\gamma_{l,s}K_0(\gamma_{l,s}r)\right)\right. \\ &\left.- \gamma_{l,p} \left(\frac{x_3^2-x_2^2}{r}K_1(\gamma_{l,p}r) - x^2\gamma_{l,p}K_0(\gamma_{l,p}r)\right) \right), \\
    c_{l,32}(x_2,x_3) &= c_{l,23}(x_2,x_3) \\&= \frac{x_2x_3}{4\pi^2\rho\omega^2r^2}\left(\frac{2}{r}\gamma_{l,p}K_1(\gamma_{l,p}r) + \gamma_{l,p}^2K_0(\gamma_{l,p}r) - \frac{2}{r}\gamma_{l,s}K_1(\gamma_{l,s}r) - \gamma_{l,s}^2K_0(\gamma_{l,s}r)\right), \\
    c_{l,33}(x_2,x_3) &= -\frac{1}{4\pi^2\mu}K_0(\gamma_{l,s}r) + \frac{1}{4\pi^2\rho\omega^2r^2}\left(\gamma_{l,s}\left(\frac{x_2^2-x_3^2}{r}K_1(\gamma_{l,s}r)-x_3^2\gamma_{l,s}K_0(\gamma_{l,s}r)\right)\right. \\
    &\left.-\gamma_{l,p}\left(\frac{x_2^2-x_3^2}{r}K_1(\gamma_{l,p}r)-x_3^2\gamma_{l,p}K_0(\gamma_{l,p}r)\right)\right).
\end{align*}
Due to the evanescent nature of $K_0(r)$ and $K_1(r)$, all $c_{l,ij}$ decay exponentially with increasing $r$. This implies that the corresponding part in Green's function generates an evanescent wave. As a result, all modes with $l$ satisfying $(\alpha+l)^2\geq k_s^2>k_p^2$ are evanescent modes.

\item \textbf{Case II: $k_s^2>(\alpha+l)^2\geq k_p^2$.}

In this case $D = \mu(\lambda+2\mu)(|\bm{\xi}|^2-\gamma_{l,s}^2)(|\bm{\xi}|^2+\gamma_{l,p}^2)$. Using Lemma \ref{lem:aux1} and Lemma \ref{lem:aux2}, we have
\begin{align*}
    c_{l,11}(x_2,x_3) &= -\frac{\ii\gamma_{l,s}^2}{8\pi\rho\omega^2}H_0^{(1)}(\gamma_{l,s}r) - \frac{(\alpha+l)^2}{4\pi^2\rho\omega^2}K_0(\gamma_{l,p}r), \\
    c_{l,21}(x_2,x_3) &= c_{l,12}(x_2,x_3)= \frac{(\alpha+l)x_2}{4\pi^2\rho\omega^2 r}\left(\frac{\pi}{2}\gamma_{l,s}H_1^{(1)}(\gamma_{l,s}r) - \ii\gamma_{l,p}K_1(\gamma_{l,p}r)\right), \\
    c_{l,31}(x_2,x_3) &= c_{l,13}(x_2,x_3)= \frac{(\alpha+l)x_3}{4\pi^2\rho\omega^2 r}\left(\frac{\pi}{2}\gamma_{l,s}H_1^{(1)}(\gamma_{l,s}r) - \ii\gamma_{l,p}K_1(\gamma_{l,p}r)\right),\\
    c_{l,22}(x_2,x_3) &= \frac{\ii}{8\pi\mu}H_0^{(1)}(\gamma_{l,s} r) - \frac{1}{4\pi^2\rho\omega^2}\left(\frac{\gamma_{l,p}}{r^2}\left(\frac{x_3^2-x_2^2}{r}K_1(\gamma_{l,p}r) - x_2^2\gamma_{l,p}K_0(\gamma_{l,p}r)\right)\right. \\ &\left.-\frac{\pi\ii \gamma_{l,s}}{2r^2}\left(\frac{x_3^2-x_2^2}{r}H_1^{(1)}(\gamma_{l,s}r) - x_2^2\gamma_{l,s}H_0^{(1)}(\gamma_{l,s}r)\right)\right), \\
    c_{l,32}(x_2,x_3) &= c_{l,23}(x_2,x_3) \\&= \frac{x_2x_3}{4\pi^2\rho\omega^2 r}\left(-\frac{1}{2}\pi\ii\gamma_{l,s}\left(-\frac{2}{r}H_1^{(1)}(\gamma_{l,s}r) + \gamma_{l,s}H_0^{(1)}(\gamma_{l,s}r)\right) + \frac{\gamma_{l,p}}{r}\left(\frac{2}{r}K_1(\gamma_{l,p}r) + \gamma_{l,p}K_0(\gamma_{l,p}r)\right)\right), \\
    c_{l,33}(x_2,x_3) &= \frac{\ii}{8\pi\mu}H_0^{(1)}(\gamma_{l,s} r) - \frac{1}{4\pi^2\rho\omega^2}\left(\frac{\gamma_{l,p}}{r^2}\left(\frac{x_2^2-x_3^2}{r}K_1(\gamma_{l,p}r) - x_3^2\gamma_{l,p}K_0(\gamma_{l,p}r)\right)\right. \\ &\left.-\frac{\pi\ii \gamma_{l,s}}{2r^2}\left(\frac{x_2^2-x_3^2}{r}H_1^{(1)}(\gamma_{l,s}r) - x_3^2\gamma_{l,s}H_0^{(1)}(\gamma_{l,s}r)\right)\right).
\end{align*}

We can see that the p-wave part in the Green's function is exponentially decreasing, while the s-wave part propagates. Therefore, when $k_s^2>(\alpha+l)^2\geq k_p^2$, only s-wave contributes to the propagating solution, and p-wave stays evanescent.

\item \textbf{Case III: $k_s^2> k_p^2 > (\alpha+l)^2$.}

In this case $D = \mu(\lambda+2\mu)(|\bm{\xi}|^2-\gamma_{l,s}^2)(|\bm{\xi}|^2-\gamma_{l,p}^2)$. Using Lemma \ref{lem:aux2}, we have
\begin{align*}
    c_{l,11}(x_2,x_3) &= -\frac{\ii}{8\pi\rho\omega^2}\gamma_{l,s}^2H_0^{(1)}(\gamma_{l,s}r) + \frac{\ii(\alpha+l)^2}{8\pi\rho\omega^2}H_0^{(1)}(\gamma_{l,p}r), \\
    c_{l,21}(x_2,x_3) &= c_{l,12}(x_2,x_3)= \frac{(\alpha+l)x_2}{8\pi\rho\omega^2 r}(\gamma_{l,s}H_1^{(1)}(\gamma_{l,s}r) - \gamma_{l,p}H_1^{(1)}(\gamma_{l,p}r)), \\
    c_{l,31}(x_2,x_3) &= c_{l,13}(x_2,x_3)= \frac{(\alpha+l)x_3}{8\pi\rho\omega^2 r}(\gamma_{l,s}H_1^{(1)}(\gamma_{l,s}r) - \gamma_{l,p}H_1^{(1)}(\gamma_{l,p}r)), \\
    c_{l,22}(x_2,x_3) &= \frac{\ii}{8\pi\mu}H_0^{(1)}(\gamma_{l,s}r) + \frac{\ii}{8\pi\rho\omega^2r^2}\left(\gamma_{l,s}\left(\frac{x_3^2-x_2^2}{r}H_1^{(1)}(\gamma_{l,s}r)-x_2^2\gamma_{l,s}H_0^{(1)}(\gamma_{l,s}r)\right)\right. \\
    &\left.-\gamma_{l,p}\left(\frac{x_3^2-x_2^2}{r}H_1^{(1)}(\gamma_{l,s}r)-x_2^2\gamma_{l,s}H_0^{(1)}(\gamma_{l,s}r)\right)\right), \\
    c_{l,32}(x_2,x_3) &= c_{l,23}(x_2,x_3) \\&= \frac{\ii x_2x_3}{8\pi\rho\omega^2r}\left(-\gamma_{l,p}\left(-\frac{2}{r}H_1^{(1)}(\gamma_{l,p}r) + \gamma_{l,p}H_0^{(1)}(\gamma_{l,p}r)\right) + \gamma_{l,s}\left(-\frac{2}{r}H_1^{(1)}(\gamma_{l,s}r) + \gamma_{l,s}H_0^{(1)}(\gamma_{l,s}r)\right)\right), \\
    c_{l,33}(x_2,x_3) &= \frac{\ii}{8\pi\mu}H_0^{(1)}(\gamma_{l,s}r) + \frac{\ii}{8\pi\rho\omega^2r^2}\left(\gamma_{l,s}\left(\frac{x_2^2-x_3^2}{r}H_1^{(1)}(\gamma_{l,s}r)-x_3^2\gamma_{l,s}H_0^{(1)}(\gamma_{l,s}r)\right)\right. \\
    &\left.-\gamma_{l,p}\left(\frac{x_2^2-x_3^2}{r}H_1^{(1)}(\gamma_{l,s}r)-x_3^2\gamma_{l,s}H_0^{(1)}(\gamma_{l,s}r)\right)\right).
\end{align*}
\end{itemize}
This leads to a propagating solution for both p-wave and s-wave instead of an evanescent one.

To conclude, we achieved the full expression of $\tgao(\bm{x})$. The desired quasi-periodic Green's function satisfying \eqref{eq:quasi-p-Green-3d} could be calculated by
\begin{align*}
&\widetilde{\bm{G}}^{\alpha,\omega}(\bm{x},\bm{y}) = e^{\ii\alpha (x_1-y_1)}\tgao(\bm{x},\bm{y}) = \sum\limits_{l\in 2\pi\mathbb{Z}}e^{\ii (\alpha+l)(x_1-y_1)}\bm{c}_l \\
&=: \sum\limits_{\substack{l\in 2\pi\mathbb{Z} \\ (\alpha+l)^2<k_p^2}} \e^{\ii(\alpha+l)(x_1-y_1)}\widetilde{\bm{G}}_1^{\alpha_l} + \sum\limits_{\substack{l\in 2\pi\mathbb{Z} \\ k_p^2\leq(\alpha+l)^2<k_s^2}} \e^{\ii(\alpha+l)(x_1-y_1)}\widetilde{\bm{G}}_2^{\alpha_l} + \sum\limits_{\substack{l\in 2\pi\mathbb{Z} \\ (\alpha+l)^2\geq k_s^2}}\e^{\ii(\alpha+l)(x_1-y_1)}\widetilde{\bm{G}}_3^{\alpha_l}.
\end{align*}
\begin{rmk}\label{rem:3d1}
From the expressions of $\widetilde{\bm{G}}^{\alpha_l}_i$, $i=1,2,3$, we could easily observe that
$$
\widetilde{\bm{G}}_i^{-\alpha_l}((y_2,y_3),(x_2,x_3)) = \widetilde{\bm{G}}_i^{\alpha_l}((x_2,x_3),(y_2,y_3)).
$$
\end{rmk}

\subsection{3D quasi-biperiodic Green's function}\label{sec:3d2}
The three-dimensional $\bm{\alpha}$-quasi-biperiodic Green's function of Lam\'{e} system  $\mathring{\bm{G}}^{\bm{\alpha},\omega}$ satisfies
\begin{equation}\label{eq:quasi-p-Green-3d2}
\begin{split}
(\mathcal{L}^{\lambda, \mu} +\rho\omega^2\bm{I}_3)\mathring{\bm{G}}^{\bm{\alpha},\omega}(\bm{x},\bm{y})&=\sum_{n_1,n_2 \in \mathbb{Z}} \delta(x_1 - y_1 - n_1) \delta(x_2 - y_2-n_2)\delta(x_3-y_3) {\rm{e}}^{{\rm{i}} \bm{n}\cdot\bm{\alpha}}\bm{I}_3.
\end{split}
\end{equation}
Here, $\bm{\alpha}:=(\alpha_1,\alpha_2)$ denotes the quasi-periodicity, $\bm{n}=(n_1,n_2)$. Since $\mathring{\bm{G}}^{\bm{\alpha},\omega}(\bm{x},\bm{y}) = \mathring{\bm{G}}^{\bm{\alpha},\omega}(\bm{x}-\bm{y})$, we regard $\mathring{\bm{G}}^{\bm{\alpha},\omega}$ as a function of $\bm{x}$ now. Recall that $\bar{\bm{x}}:=(x_1,x_2)$. Define $\doverline{\bm{G}}^{\bm{\alpha},\omega}(\bm{x}):= e^{-\ii \bm{\alpha}\cdot\bar{\bm{x}}}\mathring{\bm{G}}^{\bm{\alpha},\omega}(\bm{x})$. Then, $\doverline{\bm{G}}^{\bm{\alpha},\omega}(\bm{x})$ is a biperiodic function of $x_1$ and $x_2$ with periodicity $(1,1)$. By straightforward calculation we get
    \begin{equation}
        \begin{aligned}
            e^{-\ii\bm{\alpha}\cdot\bar{\bm{x}}}(\llm(e^{\ii\bm{\alpha}\cdot\bar{\bm{x}}} \doverline{\bm{G}}^{\bm{\alpha},\omega}(\bm{x})))_{1j} &= \mu(-|\bm{\alpha}|^2+2\ii\bm{\alpha}\cdot\bar{\nabla} + \Delta)\doverline{G}^{\bm{\alpha},\omega}_{1j}+(\lambda+\mu)(\ii\alpha_1+\partial_1)^2\doverline{G}^{\bm{\alpha},\omega}_{1j} \\ &+ (\lambda+\mu)(\ii\alpha_1+\partial_1)(\ii\alpha_2+\partial_2)\doverline{G}^{\bm{\alpha},\omega}_{2j} + (\lambda+\mu)(\ii\alpha_1\partial_3 + \partial_{13})\doverline{G}^{\bm{\alpha},\omega}_{3j}, \\
            e^{-\ii\bm{\alpha}\cdot\bar{\bm{x}}}(\llm(e^{\ii\bm{\alpha}\cdot\bar{\bm{x}}} \doverline{\bm{G}}^{\bm{\alpha},\omega}(\bm{x})))_{2j} &= \mu(-|\bm{\alpha}|^2+2\ii\bm{\alpha}\cdot\bar{\nabla} + \Delta)\doverline{G}^{\bm{\alpha},\omega}_{2j}+(\lambda+\mu)(\ii\alpha_1+\partial_1)(\ii\alpha_2+\partial_2)\doverline{G}^{\bm{\alpha},\omega}_{1j} \\ &+(\lambda+\mu)(\ii\alpha_2+\partial_2)^2\doverline{G}^{\bm{\alpha},\omega}_{2j}+(\lambda+\mu)(\ii\alpha_2\partial_3 + \partial_{23})\doverline{G}^{\bm{\alpha},\omega}_{3j}, \\
             e^{-\ii\bm{\alpha}\cdot\bar{\bm{x}}}(\llm(e^{\ii\bm{\alpha}\cdot\bar{\bm{x}}} \doverline{\bm{G}}^{\bm{\alpha},\omega}(\bm{x})))_{3j} &= \mu(-|\bm{\alpha}|^2+2\ii\bm{\alpha}\cdot\bar{\nabla} + \Delta)\doverline{G}^{\bm{\alpha},\omega}_{3j}+(\lambda+\mu)(\ii\alpha_1\partial_3+\partial_{13})\doverline{G}^{\bm{\alpha},\omega}_{1j} \\&+ (\lambda+\mu)(\ii\alpha_2\partial_3+\partial_{13})\doverline{G}^{\bm{\alpha},\omega}_{2j} + (\lambda+\mu)\partial_{33}\doverline{G}^{\bm{\alpha},\omega}_{3j}.
        \end{aligned}
    \end{equation}
    Here $j=1,2,3, k=2,3$, $\partial_if$ is short for $\frac{\partial f}{\partial x_i}$, and $\partial_{il}f$ is short for $\frac{\partial^2 f}{\partial x_i \partial x_l}$. Define $\doverline{\mathcal{L}}^{\lambda,\mu} := (d\overline{\mathcal{L}}^{\lambda,\mu}_{ij})_{3\times 3}, i,j=1,2,3$ as
    \begin{align*}
    \doverline{\mathcal{L}}^{\lambda,\mu}_{11} &= (\lambda+2\mu)(\ii\alpha_1+\partial_1)^2 + \mu(\ii\alpha_2+\partial_2)^2 + \mu\partial_{33}, \\
    \doverline{\mathcal{L}}^{\lambda,\mu}_{22} &= (\lambda+2\mu)(\ii\alpha_2+\partial_2)^2 + \mu(\ii\alpha_1+\partial_1)^2 + \mu\partial_{33}, \\
    \doverline{\mathcal{L}}^{\lambda,\mu}_{33} &= \mu(-|\bm{\alpha}|^2+2\ii\bm{\alpha}\cdot\bar{\nabla}+\Delta) + (\lambda+\mu)\partial_{33}, \\
    \doverline{\mathcal{L}}^{\lambda,\mu}_{12} &= \doverline{\mathcal{L}}^{\lambda,\mu}_{21} = (\lambda+\mu)(\ii\alpha_1+\partial_1)(\ii\alpha_2+\partial_2), \\
    \doverline{\mathcal{L}}^{\lambda,\mu}_{13} &= \doverline{\mathcal{L}}^{\lambda,\mu}_{31} = (\lambda+\mu)(\ii\alpha_1\partial_3 + \partial_{13}), \\
    \doverline{\mathcal{L}}^{\lambda,\mu}_{23} &= \doverline{\mathcal{L}}^{\lambda,\mu}_{32} = (\lambda+\mu)(\ii\alpha_2\partial_3 + \partial_{23}).
    \end{align*}
Then
\begin{equation*}
    (\doverline{\mathcal{L}}^{\lambda,\mu}+\rho\omega^2\bm{I}))\doverline{\bm{G}}^{\bm{\alpha},\omega}(\bm{x}) = \sum\limits_{n_1, n_2\in\mathbb{Z}}\delta(x_1-n_1)\delta(x_2-n_2)\delta(x_3)\bm{I}_3.
\end{equation*}
Consider the Fourier expansion of $\doverline{\bm{G}}^{\bm{\alpha},\omega}(\bm{x})$ with respect to $x_1$ and $x_2$ in the form of
\begin{equation*}
    \doverline{\bm{G}}^{\bm{\alpha},\omega}_{ij}(x_1,x_2,x_3) = \sum\limits_{l_1, l_2\in 2\pi\mathbb{Z}} c_{l,ij}(x_3)e^{\ii (l_1x_1+l_2x_2)}
\end{equation*}
and use Poisson summation formula, we get:
\begin{lem}
    $\bm{c}_l$ satisfies
    \begin{equation*}
        \begin{aligned}
            &(-(\lambda+2\mu)(\alpha_1+l_1)^2-\mu(\alpha_2+l_2)^2+\mu\partial_{33}+\rho\omega^2)c_{l,1j} -(\lambda+\mu)(\alpha_1+l_1)(\alpha_2+l_2)c_{l,1j} \\ &+ \ii(\lambda+\mu)(\alpha_1+l_1)\partial_3 c_{l,3j} = \frac{1}{4\pi^2}\delta_0\delta_{1j}, \\
            &-(\lambda+\mu)(\alpha_1+l_1)(\alpha_2+l_2)c_{l,1j} + (-(\lambda+2\mu)(\alpha_2+l_2)^2-\mu(\alpha_1+l_1)^2+\mu\partial_{33}+\rho\omega^2)c_{l,2j} \\ &+ \ii(\lambda+\mu)(\alpha+l_2)\partial_3 c_{l,3j} = \frac{1}{4\pi^2}\delta_0\delta_{2j}, \\
            &(-\mu(\alpha_1+l_1)^2-(\alpha_2+l_2)^2 +(\lambda+2\mu)\partial_{33}+\rho\omega^2)c_{l,3j}+\ii(\lambda+\mu)(\alpha_1+l_1)\partial_3c_{l,1j} \\& + \ii(\lambda+\mu)(\alpha_2+l_2)\partial_3c_{l,2j}  = \frac{1}{4\pi^2}\delta_0\delta_{3j}
        \end{aligned}
    \end{equation*}
    for $j=1,2,3$.
\end{lem}
Define $\bm{\alpha}_l:=(\alpha_1+l_1, \alpha_2+l_2)$, $\alpha_{l,1}=\alpha_1+l_1$ and $\alpha_{l,2}=\alpha_2+l_2$. By applying Fourier transform, we get
\begin{equation*}
\left\{
\begin{aligned}
(-(\lambda+2\mu)\alpha_{l,1}^2-\mu\alpha_{l,2}^2-\mu\xi^2+\rho\omega^2)\hat{c}_{l,1j} - (\lambda+\mu)\alpha_{l,1}\alpha_{l,2}\hat{c}_{l,2j} - (\lambda+\mu)\alpha_{l,1}\xi \hat{c}_{l,3j}=(2\pi)^{-5/2}\delta_{1j}, \\
-(\lambda+\mu)\alpha_{l,1}\alpha_{l,2}\hat{c}_{l,1j} + (-(\lambda+2\mu)\alpha_{l,2}^2-\mu\alpha_{l,1}^2-\mu\xi^2+\rho\omega^2)\hat{c}_{l,2j}-(\lambda+\mu)\alpha_{l,2}\xi \hat{c}_{l,3j} = (2\pi)^{-5/2}\delta_{2j}, \\
-(\lambda+\mu)\alpha_{l,1}\xi \hat{c}_{l,1j} - (\lambda+\mu)\alpha_{l,2}\xi \hat{c}_{l,2j} + (-\mu(\alpha_{l,1}^2 + \alpha_{l,2}^2)-(\lambda+2\mu)\xi^2+\rho\omega^2)\hat{c}_{l,3j} = (2\pi)^{-5/2}\delta_{3j}.
\end{aligned}
\right.
\end{equation*}
Define $\kappa^2 := \alpha_{l,1}^2+\alpha_{l,2}^2+\xi^2 =: |\bm{\alpha}_l|^2+\xi^2$, $\gamma_{l,p} = \sqrt{||\bm{\alpha}_l|^2-\frac{\rho\omega^2}{\lambda+2\mu}|}$, $\gamma_{l,s} = \sqrt{||\bm{\alpha}_l|^2 - \frac{\rho\omega^2}{\mu}|}$. The straightforward calculation yields
\begin{align*}
    \hat{c}_{l,11}(\xi) &= \frac{1}{(2\pi)^{5/2}\kappa^2}\left(\frac{\alpha_{l,2}^2+\xi^2}{\rho\omega^2-\mu\kappa^2} + \frac{\alpha_{l,1}^2}{\rho\omega^2-(\lambda+2\mu)\kappa^2}\right), \\
    \hat{c}_{l,21}(\xi) &= \frac{\alpha_{l,1}\alpha_{l,2}}{(2\pi)^{5/2}\kappa^2}\left(\frac{1}{\rho\omega^2-(\lambda+2\mu)\kappa^2} - \frac{1}{\rho\omega^2-\mu\kappa^2}\right) = \hat{c}_{l,12}, \\
    \hat{c}_{l,31}(\xi) &= \frac{\alpha_{l,1}\xi}{(2\pi)^{5/2}\kappa^2}\left(\frac{1}{\rho\omega^2-(\lambda+2\mu)\kappa^2}-\frac{1}{\rho\omega^2-\mu\kappa^2}\right) = \hat{c}_{l,13}, \\
    \hat{c}_{l,22}(\xi) &= \frac{1}{(2\pi)^{5/2}\kappa^2} \left(\frac{\alpha_{l,2}^2}{ \rho\omega^2 - (\lambda+2\mu)\kappa^2} + \dfrac{\alpha_{l,1}^2 + \xi^2}{\rho\omega^2 - \mu\kappa^2} \right), \\
    \hat{c}_{l,32}(\xi) &= \frac{\alpha_{l,2} \xi}{(2\pi)^{5/2}\kappa^2} \left( \frac{1}{\rho\omega^2 - (\lambda+2\mu)\kappa^2} - \frac{1}{\rho\omega^2 - \mu\kappa^2} \right) = \hat{c}_{l,23}, \\
    \hat{c}_{l,33}(\xi) &= \frac{1}{(2\pi)^{5/2}\kappa^2} \left( \frac{\xi^2}{\rho\omega^2 - (\lambda+2\mu)\kappa^2} + \frac{\alpha_{l,1}^2 + \alpha_{l,2}^2}{\rho\omega^2 - \mu\kappa^2} \right).
\end{align*}

Using some known results on Inverse Fourier transform (cf. Lemma 2 and Lemma 3 in \cite{WuHe}, we can calculate explicit expressions of $c_{l,ij}$ in three different cases.
\begin{itemize}
\item \textbf{Case I: $|\bm{\alpha}_l|^2\geq\rho\omega^2/\mu>\rho\omega^2/(\lambda+2\mu)$.}
In this case, we have
\begin{align*}
    c_{l,ii} &= \frac{1}{8\pi^2}\left(-\frac{1}{\mu}\frac{\e^{-\gamma_{l,s}|x_3|}}{\gamma_{l,s}} + \frac{\alpha_{l,i}^2}{\rho\omega^2}\left(\frac{\e^{-\gamma_{l,s}|x_3|}}{\gamma_{l,s}} - \frac{\e^{-\gamma_{l,p}|x_3|}}{\gamma_{l,p}}\right)\right), ~i=1,2, \\
    c_{l,12} &= \frac{\alpha_{l,1}\alpha_{l,2}}{8\pi^2\rho\omega^2}\left(\frac{\e^{-\gamma_{l,s}|x_3|}}{\gamma_{l,s}} - \frac{\e^{-\gamma_{l,p}|x_3|}}{\gamma_{l,p}}\right) = c_{l,21}, \\
    c_{l,i3} &= \frac{\ii\alpha_{l,i}}{8\pi^2\rho\omega^2}\mathrm{sgn}(x_3)(\e^{-\gamma_{l,s}|x_3|} - \e^{-\gamma_{l,p}|x_3|}) = c_{l,3i},~i=1,2,\\
    c_{l,33} &= \frac{1}{8\pi^2}\left(-\frac{1}{\mu}\frac{\e^{-\gamma_{l,s}|x_3|}}{\gamma_{l,s}} + \frac{1}{\rho\omega^2}\left(\gamma_{l,p}\e^{-\gamma_{l,p}|x_3|}-\gamma_{l,s}\e^{-\gamma_{l,s}|x_3|}\right)\right).
\end{align*}
The result implies that, both s-wave part and p-wave part are evanescent waves in $x_3$ direction.
\item \textbf{Case II: $\rho\omega^2/\mu>|\bm{\alpha}_l|^2\geq \rho\omega^2/(\lambda+2\mu)$.}
In this case, we have
\begin{align*}
    c_{l,ii} &= \frac{1}{8\pi^2}\left(-\frac{1}{\mu}\frac{\ii\e^{\ii\gamma_{l,s}|x_3|}}{\gamma_{l,s}} + \frac{\alpha_{l,i}^2}{\rho\omega^2}\left(\frac{\ii\e^{\ii\gamma_{l,s}|x_3|}}{\gamma_{l,s}} - \frac{\e^{-\gamma_{l,p}|x_3|}}{\gamma_{l,p}}\right)\right), ~i=1,2, \\
    c_{l,12} &= \frac{\alpha_{l,1}\alpha_{l,2}}{8\pi^2\rho\omega^2}\left(\frac{\ii\e^{\ii\gamma_{l,s}|x_3|}}{\gamma_{l,s}} - \frac{\e^{-\gamma_{l,p}|x_3|}}{\gamma_{l,p}}\right) = c_{l,21}, \\
    c_{l,i3} &= \frac{\ii\alpha_{l,i}}{8\pi^2\rho\omega^2}\mathrm{sgn}(x_3)(\e^{\ii\gamma_{l,s}|x_3|} - \e^{-\gamma_{l,p}|x_3|}) = c_{l,3i},~i=1,2,\\
    c_{l,33} &= \frac{1}{8\pi^2}\left(-\frac{1}{\mu}\frac{\ii\e^{\ii\gamma_{l,s}|x_3|}}{\gamma_{l,s}} + \frac{1}{\rho\omega^2}\left(\gamma_{l,p}\e^{-\gamma_{l,p}|x_3|}-\ii\gamma_{l,s}\e^{\ii\gamma_{l,s}|x_3|}\right)\right).
\end{align*}
The result implies that, s-wave part becomes propagating wave, while p-wave part remains evanescent in $x_3$ direction.
\item \textbf{Case III: $|\bm{\alpha}_l|^2<\rho\omega^2/(\lambda+2\mu)$.}
In this case, we have
\begin{align*}
    c_{l,ii} &= \frac{1}{8\pi^2}\left(-\frac{1}{\mu}\frac{\ii\e^{\ii\gamma_{l,s}|x_3|}}{\gamma_{l,s}} + \frac{\alpha_{l,i}^2}{\rho\omega^2}\left(\frac{\ii\e^{\ii\gamma_{l,s}|x_3|}}{\gamma_{l,s}} - \frac{\ii\e^{\ii\gamma_{l,p}|x_3|}}{\gamma_{l,p}}\right)\right), ~i=1,2, \\
    c_{l,12} &= \frac{\ii\alpha_{l,1}\alpha_{l,2}}{8\pi^2\rho\omega^2}\left(\frac{\e^{\ii\gamma_{l,s}|x_3|}}{\gamma_{l,s}} - \frac{\e^{\ii\gamma_{l,p}|x_3|}}{\gamma_{l,p}}\right) = c_{l,21}, \\
    c_{l,i3} &= \frac{\ii\alpha_{l,i}}{8\pi^2\rho\omega^2}\mathrm{sgn}(x_3)(\e^{\ii\gamma_{l,s}|x_3|} - \e^{\ii\gamma_{l,p}|x_3|}) = c_{l,3i},~i=1,2,\\
    c_{l,33} &= \frac{\ii}{8\pi^2}\left(-\frac{1}{\mu}\frac{\e^{\ii\gamma_{l,s}|x_3|}}{\gamma_{l,s}} + \frac{1}{\rho\omega^2}\left(\gamma_{l,p}\e^{\ii\gamma_{l,p}|x_3|}-\gamma_{l,s}\e^{\ii\gamma_{l,s}|x_3|}\right)\right).
\end{align*}
The result implies that both s-wave part and p-wave part are propagating wave in $x_3$ direction. 
\end{itemize}

Till now we have achieved the full expression of $\doverline{\bm{G}}^{\bm{\alpha},\omega}(\bm{x})$. The desired quasi-biperiodic Green's function satisfying \eqref{eq:quasi-p-Green-3d2} could be calculated by
\begin{align*}
&\mathring{\bm{G}}^{\bm{\alpha},\omega}(\bm{x},\bm{y}) = \e^{\ii\overline{\bm{x}-\bm{y}}\cdot\bm{\alpha}}\doverline{\bm{G}}^{\bm{\alpha},\omega}(\bm{x},\bm{y}) = \sum\limits_{l_1,l_2\in 2\pi\mathbb{Z}}\e^{\ii (\alpha_1+l_1)(x_1-y_1)}\e^{\ii(\alpha_2+l_2)(x_2-y_2)}\bm{c}_l \\
&=: \sum\limits_{|\bm{\alpha}_l|^2<k_p^2} \e^{\ii\alpha_l\cdot\overline{\bm{x}-\bm{y}}}\mathring{\bm{G}}_1^{\alpha_l} + \sum\limits_{k_p^2\leq|\bm{\alpha}_l|^2<k_s^2} \e^{\ii\alpha_l\cdot\overline{\bm{x}-\bm{y}}}\mathring{\bm{G}}_2^{\alpha_l} + \sum\limits_{|\bm{\alpha}_l|^2>k_s^2}\e^{\ii\alpha_l\cdot\overline{\bm{x}-\bm{y}}}\mathring{\bm{G}}_3^{\alpha_l}.
\end{align*}
\begin{rmk}\label{rem:3d2}
From the expressions of $\mathring{\bm{G}}^{\alpha_l}_i$, $i=1,2,3$, we could easily observe that
$$
\mathring{\bm{G}}_i^{-\alpha_l}(\bm{y},\bm{x}) = \mathring{\bm{G}}_i^{\alpha_l}(\bm{x},\bm{y}).
$$
\end{rmk}


\section{General properties of quasi-periodic/biperiodic Green's functions}\label{sec4}
In this section, we give reciprocity relations for point sources, scattered fields, and total fields, respectively. Since the frequency $\omega$ stays the same in our discussion, for brevity we will leave out $\omega$ when mentioning Green's functions. Let $\hat{\alpha}=-\alpha$ and $\hat{\alpha}_l=l+\hat{\alpha}$. Then $\hat{\alpha}_{-l}=-\alpha_l$, $\hat{\alpha}^2_{-l}=\alpha^2_l$. We can directly get the following reciprocity result.

\begin{lem}[Reciprocity relation for point source]\label{lem:31}
$$
\bm{G}^{\alpha}(\bm{x},\bm{z})=\bm{G}^{\hat{\alpha}}(\bm{z},\bm{x}),\ \bm{x}\neq \bm{z},\ \bm{x},\bm{z}\in\mathbb{R}^2.
$$
The same results hold for $\widetilde{\bm{G}}^{\alpha}$ and $\mathring{\bm{G}}^{\alpha}$.
\end{lem}

\begin{proof}
From Section \ref{sec:green2d} and Remark \ref{rem:2d}, we have

\begin{equation}\label{eq:greenfunction}
    \bm{G}^{\hat{\alpha}}(\bm{z},\bm{x}) = \sum\limits_{\substack{l\in 2\pi\mathbb{Z} \\ \hat{\alpha}_l^2<k_p^2}} \e^{\ii\hat{\alpha}_l(z_1-x_1)}\bm{G}_1^{\hat{\alpha}_l} + \sum\limits_{\substack{l\in 2\pi\mathbb{Z} \\ k_p^2\leq\hat{\alpha}_l^2<k_s^2}} \e^{\ii\hat{\alpha}_l(z_1-x_1)}\bm{G}_2^{\hat{\alpha}_l} + \sum\limits_{\substack{l\in 2\pi\mathbb{Z} \\ \hat{\alpha}_l^2\geq k_s^2}} \e^{\ii\hat{\alpha}_l(z_1-x_1)}\bm{G}_3^{\hat{\alpha}_l}.
\end{equation}

Notice that
\begin{align*}
    &\sum\limits_{\substack{l\in 2\pi\mathbb{Z} \\ \hat{\alpha}_l^2<k_p^2}}\e^{\ii\hat{\alpha}_l(z_1-x_1)}\bm{G}_1^{\hat{\alpha}_l}(z_2-x_2) = \sum\limits_{\substack{l\in 2\pi\mathbb{Z} \\ \alpha_{-l}^2<k_p^2}} \e^{\ii\alpha_{-l}(x_1-z_1)}\bm{G}_1^{-\alpha_{-l}}(z_2-x_2) \\ &= \sum\limits_{\substack{l\in 2\pi\mathbb{Z} \\ \alpha_{-l}^2<k_p^2}}\e^{\ii\alpha_{-l}(x_1-z_1)}\bm{G}_1^{\alpha_{-l}}(x_2-z_2) = \sum\limits_{\substack{l\in 2\pi\mathbb{Z} \\ \alpha_{l}^2<k_p^2}}\e^{\ii\alpha_{l}(x_1-z_1)}\bm{G}_1^{\alpha_{l}}(x_2-z_2).
\end{align*}
The last equality holds by replacing the summation dummy variable with $-l$ in the previous summation. By applying the same method to the other two terms in \eqref{eq:greenfunction}, we immediately get the desired result. The results for $\widetilde{\bm{G}}^\alpha$ and $\mathring{\bm{G}}^\alpha$ follow from Section \ref{sec:3d1}, Remark \ref{rem:3d1} and Section \ref{sec:3d2}, Remark \ref{rem:3d2}.
\end{proof}

\begin{lem}[Reciprocity relation for scattered field]\label{lem:32}
For the periodic structure scattering problem \eqref{eq:Navier}-\eqref{eq:boundaryC} with point source incidence, the scattered field satisfies the following reciprocity relation
\[
\bm{p}\cdot\bm{u}_{\alpha, \mathrm{sc}}(\bm{x},\bm{z},\bm{q})=\bm{q}\cdot\bm{u}_{\hat{\alpha}, \mathrm{sc}}(\bm{z},\bm{x},\bm{p})
\]
for $\forall \bm{x},\bm{z}\in D_1$, $\bm{x}\neq\bm{z}$. The same result also holds for the quasi-periodic/biperiodic system in $\mathbb{R}^3$.
\end{lem}

\begin{proof}
Here we prove the lemma in $\mathbb{R}^2$. The proof could also be applied to three dimensional case, except a few places which we will point out explicitly. Applying the third Betti formula to $\bm{u}_{\alpha, \mathrm{sc}}(\bm{x},\bm{y},\bm{q})$ and $\bm{G}^{\alpha}(\bm{x},\bm{y})\bm{q}$ in $D_1$, we have
\begin{equation}\label{eq:L1}
\bm{u}_{\alpha, \mathrm{sc}}(\bm{x},\bm{z},\bm{q}) = \int_{\partial D_1}\big\{[T_{\bm{\nu}(\bm{y})} \bm{G}^{\alpha}(\bm{x},\bm{y})]^{\top} \bm{u}_{\alpha, \mathrm{sc}}(\bm{y},\bm{z},\bm{q}) - \bm{G}^{\alpha}(\bm{x},\bm{y})[T_{\bm{\nu}(\bm{y})} \bm{u}_{\alpha, \mathrm{sc}}(\bm{y},\bm{z},\bm{q})] \big\}\mathrm{d}s(\bm{y})
\end{equation}
for $\bm{x}\in D_1$. Here, $T_{\bm{\nu}}$ is the surface traction operator defined by $T_{\bm{\nu}}(\bm{u}):=2\mu {\partial}_{\bm{\nu}} \bm{u}+\lambda \bm{\nu} \nabla \cdot \bm{u}-\mu \bm{\tau} \mathrm{curl}\bm{u}$. Due to the quasi-periodicity and the opposite directions of the normal vectors on $\Gamma_{L,1}$ and $\Gamma_{R,1}$, the boundary integral of \eqref{eq:L1} on $\Gamma_{L,1}$ and $\Gamma_{R,1}$ cancel each other.


The Green's function $\bm{G}^\alpha$(as well as $\widetilde{\bm{G}}^\alpha$ and $\mathring{\bm{G}}^{\bm{\alpha}}$) could all be written in the form of $\sum_{l\in 2\pi\mathbb{Z},i}\e^{\ii \alpha_l(x_1-y_1)}\bm{G}^{\alpha_l}_i$, where $\bm{G}^{\alpha_l}_i$ is independent of $x_1$. Since the normal vector of $\Gamma_{h,1}\cap\partial D_1$ is $(0,1)^\top$ in 2-dimensional case ($(0,0,1)^\top$ in 3-dimensional case), we observe that $x_1$-component is not involved in $T_{\bm{\nu}(\bm{y})}\bm{G}^\alpha$. Therefore, $T_{\bm{\nu}(\bm{y})}\bm{G}^\alpha$ could still be written in the form of $\sum_{l\in 2\pi\mathbb{Z},i}\e^{\ii \alpha_l(x_1-y_1)}{\bm{G}^{\alpha_l}_i}'$. Since $\bm{u}_{\alpha,\mathrm{sc}}$ satisfies the Rayleigh expansion \eqref{eq:sc_radiation}(\eqref{eq:rayleigh31} and \eqref{eq:rayleigh32} for $\mathbb{R}^3$ cases), by the orthogonality of the exponential functions over the integral of $x_1$ on $\Gamma_{h,1}\cap \partial D_1$ we conclude that
\begin{equation}
\int_{\Gamma_{h,1}\cap\partial D_1} \big\{[T_{\bm{\nu}(\bm{y})} \bm{G}^{\alpha}(\bm{x},\bm{y})]^{\top} \bm{u}_{\alpha, \mathrm{sc}}(\bm{y},\bm{z},\bm{q}) - \bm{G}^{\alpha}(\bm{x},\bm{y})[T_{\bm{\nu}(\bm{y})} \bm{u}_{\alpha, \mathrm{sc}}(\bm{y},\bm{z},\bm{q})] \big\}\mathrm{d}s(\bm{y})=0.
\end{equation}
Therefore, \eqref{eq:L1} becomes
\begin{equation}\label{eq:L2}
\bm{u}_{\alpha, \mathrm{sc}}(\bm{x},\bm{z},\bm{q}) = \int_{ \Gamma_1}\big\{[T_{\bm{\nu}(\bm{y})} \bm{G}^{\alpha}(\bm{x},\bm{y})]^{\top} \bm{u}_{\alpha, \mathrm{sc}}(\bm{y},\bm{z},\bm{q}) - \bm{G}^{\alpha}(\bm{x},\bm{y})[T_{\bm{\nu}(\bm{y})} \bm{u}_{\alpha, \mathrm{sc}}(\bm{y},\bm{z},\bm{q})] \big\}\mathrm{d}s(\bm{y}).
\end{equation}
As both the incident field and scattered field can be decomposed of compressional part and shear part (cf. \cite{ChengDong}) and using \eqref{eq:L2} and Lemma~\ref{lem:31}, we have
\begin{equation}\label{eq:T1}
\bm{p}\cdot \bm{u}_{\alpha, \mathrm{sc}}(\bm{x},\bm{z},\bm{q}) = \int_{\Gamma_1}\big\{ [T_{\bm{\nu}(\bm{y})} \bm{G}^{\alpha}(\bm{x},\bm{y})\bm{p}]\cdot \bm{u}_{\alpha, \mathrm{sc}}(\bm{y},\bm{z},\bm{q}) - \bm{G}^{\alpha}(\bm{x},\bm{y})\bm{p}\cdot [T_{\bm{\nu}(\bm{y})} \bm{u}_{\alpha, \mathrm{sc}}(\bm{y},\bm{z},\bm{q}) ] \big\}\mathrm{d}s(\bm{y})
\end{equation}
and
\begin{align}\label{eq:T2}
\nonumber \bm{q}\cdot \bm{u}_{\hat{\alpha}, \mathrm{sc}}(\bm{z},\bm{x},\bm{p}) =& \int_{\Gamma_1}\big\{ [T_{\bm{\nu}(\bm{y})} \bm{G}^{\hat{\alpha}}(\bm{z},\bm{y})\bm{q}]\cdot \bm{u}_{\hat{\alpha}, \mathrm{sc}}(\bm{y},\bm{x},\bm{p}) - \bm{G}^{\hat{\alpha}}(\bm{z},\bm{y})\bm{q}\cdot [T_{\bm{\nu}(\bm{y})} \bm{u}_{\hat{\alpha}, \mathrm{sc}}(\bm{y},\bm{x},\bm{p}) ] \big\}\mathrm{d}s(\bm{y})\\
=& \int_{\Gamma_1}\big\{ [T_{\bm{\nu}(\bm{y})} \bm{G}^{\alpha}(\bm{y},\bm{z})\bm{q}]\cdot \bm{u}_{\hat{\alpha}, \mathrm{sc}}(\bm{y},\bm{x},\bm{p}) - \bm{G}^{{\alpha}}(\bm{y},\bm{z})\bm{q}\cdot [T_{\bm{\nu}(\bm{y})} \bm{u}_{\hat{\alpha}, \mathrm{sc}}(\bm{y},\bm{x},\bm{p}) ] \big\}\mathrm{d}s(\bm{y}).
\end{align}
Since both $\bm{u}_{\alpha, \mathrm{sc}}(\bm{y},\bm{z},\bm{q})$ and $\bm{u}_{\hat{\alpha}, \mathrm{sc}}(\bm{y},\bm{x},\bm{p})$ are radiating solution, we have
\begin{align*}
    0 &= \int_{D_1} \mathcal{L}^{\lambda,\mu}_y\bm{u}_{\alpha,\mathrm{sc}}(\bm{y},\bm{z},\bm{q})\cdot\bm{u}_{\hat{\alpha},\mathrm{sc}}(\bm{y},\bm{x},\bm{p}) - \mathcal{L}^{\lambda,\mu}_y\bm{u}_{\hat{\alpha},\mathrm{sc}}(\bm{y},\bm{x},\bm{p})\cdot\bm{u}_{\alpha,\mathrm{sc}}(\bm{y},\bm{z},\bm{q})\mathrm{d}s(\bm{y}) \\
    &= \int_{\Gamma_{h,1}\cup\Gamma_{L,1}\cup\Gamma_1\cup\Gamma_{R,1}}[T_{\bm{\nu}(\bm{y})} \bm{u}_{\alpha, \mathrm{sc}}(\bm{y},\bm{z},\bm{q}) ]\cdot \bm{u}_{\hat{\alpha}, \mathrm{sc}}(\bm{y},\bm{x},\bm{p}) - [T_{\bm{\nu}(\bm{y})} \bm{u}_{\hat{\alpha}, \mathrm{sc}}(\bm{y},\bm{x},\bm{p})]\cdot \bm{u}_{\alpha, \mathrm{sc}}(\bm{y},\bm{z},\bm{q}) \mathrm{d}s(\bm{y}).
\end{align*}
Similarly as before, by the orthogonality relation of exponential functions, the above integration on $\Gamma_h$ is 0, and by the quasi-periodicity, the above integration on $\Gamma_{L,1}$ and $\Gamma_{R,1}$ cancel each other. We then conclude that
\begin{align}\label{eq:T3}
\int_{\Gamma_1}\big\{[T_{\bm{\nu}(\bm{y})} \bm{u}_{\alpha, \mathrm{sc}}(\bm{y},\bm{z},\bm{q}) ]\cdot \bm{u}_{\hat{\alpha}, \mathrm{sc}}(\bm{y},\bm{x},\bm{p}) - [T_{\bm{\nu}(\bm{y})} \bm{u}_{\hat{\alpha}, \mathrm{sc}}(\bm{y},\bm{x},\bm{p})]\cdot \bm{u}_{\alpha, \mathrm{sc}}(\bm{y},\bm{z},\bm{q}) \big\} \mathrm{d}s(\bm{y})=0.
\end{align}
Now applying Betti's formula to $\bm{G}^{\alpha}(\bm{x},\bm{y})\bm{p}$ and $\bm{G}^{\alpha}(\bm{y},\bm{z})\bm{q}$ in
\[
\Omega_{\tilde{h},1}:=\{(x_1,x_2)\in\mathbb{R}^2: 0\leq x_1\leq 2\pi, \tilde{h}<x_2<f(x_1), \tilde{h}<\min_{x_1\in [0,2\pi]} f(x_1)\}
\]
and repeating the same discussions above yields
\begin{align}\label{eq:T4}
\int_{\Gamma_1}\big\{[T_{\bm{\nu}(\bm{y})} \bm{G}^{\alpha}(\bm{x},\bm{y})\bm{p}]\cdot \bm{G}^{\alpha}(\bm{y},\bm{z})\bm{q} -[T_{\bm{\nu}(\bm{y})} \bm{G}^{\alpha}(\bm{y},\bm{z})\bm{q}]\cdot \bm{G}^{\alpha}(\bm{x},\bm{y})\bm{p}  \big\}\mathrm{d}s(\bm{y})=0.
\end{align}
Combining \eqref{eq:T1}-\eqref{eq:T4}, we have
\begin{align*}
\bm{p}\cdot &\bm{u}_{\alpha, \mathrm{sc}}(\bm{x},\bm{z},\bm{q})- \bm{q}\cdot \bm{u}_{\hat{\alpha}, \mathrm{sc}}(\bm{z},\bm{x},\bm{p}) \\
=&\int_{\Gamma_1}\big\{[T_{\bm{\nu}(\bm{y})} \bm{G}^{\alpha}(\bm{x},\bm{y})\bm{p}]\cdot \bm{u}_{\alpha, \mathrm{sc}}(\bm{y},\bm{z},\bm{q}) - \bm{G}^{\alpha}(\bm{x},\bm{y})\bm{p}\cdot [T_{\bm{\nu}(\bm{y})} \bm{u}_{\alpha, \mathrm{sc}}(\bm{y},\bm{z},\bm{q}) ]\\
& - [T_{\bm{\nu}(\bm{y})} \bm{G}^{\alpha}(\bm{y},\bm{z})\bm{q} ]\cdot \bm{u}_{\hat{\alpha}, \mathrm{sc}}(\bm{y},\bm{x},\bm{p}) + \bm{G}^{\alpha}(\bm{y},\bm{z})\bm{q}\cdot [T_{\bm{\nu}(\bm{y})}\bm{u}_{\hat{\alpha}, \mathrm{sc}}(\bm{y},\bm{x},\bm{p}) ] \big\}\mathrm{d}s(\bm{y})\\
=&\int_{\Gamma_1}\big\{[T_{\bm{\nu}(\bm{y})} \bm{G}^{\alpha}(\bm{x},\bm{y})\bm{p} ]\cdot ( \bm{u}_{\alpha, \mathrm{tot}}(\bm{y},\bm{z},\bm{q}) - \bm{G}^{\alpha}(\bm{y},\bm{z})\bm{q} ) \\&- \bm{G}^{\alpha}(\bm{x},\bm{y})\bm{p}\cdot (T_{\bm{\nu}(\bm{y})}\bm{u}_{\alpha, \mathrm{tot}}(\bm{y},\bm{z},\bm{q}) - T_{\bm{\nu}(\bm{y})} \bm{G}^{\alpha}(\bm{y},\bm{z})\bm{q} )\\
& - (T_{\bm{\nu}(\bm{y})}\bm{u}_{\alpha, \mathrm{tot}}(\bm{y},\bm{z},\bm{q})-T_{\bm{\nu}(\bm{y})}\bm{u}_{\alpha, \mathrm{sc}}(\bm{y},\bm{z},\bm{q}) )\cdot \bm{u}_{\hat{\alpha}, \mathrm{sc}}(\bm{y},\bm{x},\bm{p}) \\
&+ (\bm{u}_{\alpha, \mathrm{tot}}(\bm{y},\bm{z},\bm{q})-\bm{u}_{\alpha, \mathrm{sc}}(\bm{y},\bm{z},\bm{q}))\cdot [T_{\bm{\nu}(\bm{y})}\bm{u}_{\hat{\alpha}, \mathrm{sc}}(\bm{y},\bm{x},\bm{p})  ] \big\}\mathrm{d}s(\bm{y})\\
=&\int_{\Gamma_1}\big\{-[T_{\bm{\nu}(\bm{y})}\bm{u}_{\alpha, \mathrm{tot}}(\bm{y},\bm{z},\bm{q})]\cdot \bm{u}_{\hat{\alpha}, \mathrm{sc}}(\bm{y},\bm{x},\bm{p}) +\bm{u}_{{\alpha}, \mathrm{tot}}(\bm{y},\bm{z},\bm{q})\cdot [T_{\bm{\nu}(\bm{y})} \bm{u}_{\hat{\alpha}, \mathrm{sc}}(\bm{y},\bm{x},\bm{p})]\\
&+[T_{\bm{\nu}(\bm{y})} \bm{G}^{\alpha}(\bm{x},\bm{y})\bm{p}]\cdot \bm{u}_{\alpha, \mathrm{tot}}(\bm{y},\bm{z},\bm{q}) -\bm{G}^{\alpha}(\bm{x},\bm{y})\bm{p}\cdot [T_{\bm{\nu}(\bm{y})}\bm{u}_{\alpha, \mathrm{tot}}(\bm{y},\bm{z},\bm{q}) ] \big\}\mathrm{d}s(\bm{y})\\
=&\int_{\Gamma_1}\big\{  \bm{u}_{\alpha, \mathrm{tot}}(\bm{y},\bm{z},\bm{q})\cdot [T_{\bm{\nu}(\bm{y})}\bm{u}_{\hat{\alpha}, \mathrm{tot}}(\bm{y},\bm{x},\bm{p})] - \bm{u}_{\hat{\alpha}, \mathrm{tot}}(\bm{y},\bm{x},\bm{p})\cdot [T_{\bm{\nu}(\bm{y})}\bm{u}_{{\alpha}, \mathrm{tot}}(\bm{y},\bm{z},\bm{q}) ] \big\}\mathrm{d}s(\bm{y}).
\end{align*}
The proof is complete by using the boundary condition.
\end{proof}

\begin{lem}[Reciprocity relation for total field]\label{lem:33}
For the periodic structure scattering problem \eqref{eq:Navier}-\eqref{eq:boundaryC} with point source incidence, the total field satisfies the following reciprocity relation
\[
\bm{p}\cdot\bm{u}_{\alpha, \mathrm{tot}}(\bm{x},\bm{z},\bm{q})=\bm{q}\cdot\bm{u}_{\hat{\alpha}, \mathrm{tot}}(\bm{z},\bm{x},\bm{p})
\]
for $\forall \bm{x},\bm{z}\in D_1$, $\bm{x}\neq\bm{z}$. The same result also holds for quasi-periodic/biperiodic system in $\mathbb{R}^3$.
\end{lem}

\begin{proof}
The proof follows from Lemmas~\ref{lem:31} and~\ref{lem:32} immediately.
\end{proof}

\section{Uniqueness theorem of quasi-periodic/biperiodic system}\label{sec5}
\subsection{Uniqueness of quasi-periodic elastic system in $\mathbb{R}^2$}

In this section, we present the uniqueness of the phaseless inverse periodic problem. Assume that $f_1, \tilde{f}_1\in C^2(\mathbb{R})$ are two 1-periodic functions and $h> \max_{x_1\in[0,1]}\{f_1(x_1), \tilde{f}_1(x_1) \}$. Let $\Gamma_{j,1}$, $\Omega_{j,1}^+$, $D_{j,1}, j=1,2$ be defined as in Section 2 with $f_1$ and $\tilde{f}_1$, respectively. We denote the total fields corresponding to the impenetrable surface $f_1$ and $\tilde{f}_1$ by $\bm{u}^1_{\alpha,{\rm{tot}}}$ and $\bm{u}^2_{\alpha,{\rm{tot}}}$, respectively. Then, we have the following uniqueness result.

\begin{thm}\label{thm:2d}
Let $f_1, \tilde{f}_1\in C^2(\mathbb{R})$ be 1-periodic functions. Suppose that the $\alpha$-quasi-periodic total fields satisfy
\begin{align}
\label{eq:uni1} |\bm{p}_k\cdot \bm{u}^1_{\alpha,{\rm{tot}}}(\bm{x},\tilde{\bm{z}},\bm{q})| =& |\bm{p}_k\cdot \bm{u}^2_{\alpha,{\rm{tot}}}(\bm{x},\tilde{\bm{z}},\bm{q})|, \quad \forall \bm{x} \in \overline{\Gamma}_{h,1},\\
\label{eq:uni2} |\bm{p}_k\cdot \bm{u}^1_{\alpha,{\rm{tot}}}(\bm{x},\bm{z},\bm{q}_l)| = & |\bm{p}_k\cdot \bm{u}^2_{\alpha,{\rm{tot}}}(\bm{x},\bm{z},\bm{q}_l)|, \quad \forall(\bm{x},\bm{z}) \in \overline{\Gamma}_{h,1}\times \Sigma_1,\\
\label{eq:uni3} |\bm{p}_k\cdot \left(\bm{u}^1_{\alpha,{\rm{tot}}}(\bm{x},\tilde{\bm{z}},\bm{q})+ \bm{u}^1_{\alpha,{\rm{tot}}}(\bm{x},\bm{z},\bm{q}_l) \right)|  =& |\bm{p}_k\cdot \left(\bm{u}^2_{\alpha,{\rm{tot}}}(\bm{x},\tilde{\bm{z}},\bm{q})+ \bm{u}^2_{\alpha,{\rm{tot}}}(\bm{x},\bm{z},\bm{q}_l) \right)|,  \ \forall(\bm{x}, \bm{z}) \in \overline{\Gamma}_{h,1}\times\Sigma_1,
\end{align}
for any fixed $\tilde{\bm{z}}\in\Gamma_{0,1}$, open subset $\overline{\Gamma}_{h,1}\subset\Gamma_{h,1}$, $\bm{q}, \bm{q}_l\in\mathbb{S}$, $l=1,2$, $k=1,2,\cdots,5$ and the wavenumbers $k_p\in [k_{\mathrm{min}}, k_{\mathrm{max}}]$ for some $0<k_{\mathrm{min}}<k_{\mathrm{max}}$. Here, $\Sigma_1$ is admissible curve with respect to $\Omega_1$, and $\bm{q}_1,\bm{q}_2$ and $\bm{p}_1,\bm{p}_2,\cdots,\bm{p}_5$ are any given polarization vectors that $\bm{q}_1, \bm{q}_2$ and $\bm{p}_i, \bm{p}_j$ are not colinear, $i\neq j$ ranging from 1 to 5. Then $f_1$ and $\tilde{f}_1$ coincide.

\end{thm}

\begin{proof}
It follows from \eqref{eq:uni1}-\eqref{eq:uni3} that
\begin{equation}\label{eq:realpart}
\mathrm{Re}\left(\bm{p}_k\cdot \bm{u}^1_{\alpha,{\rm{tot}}}(\bm{x},\tilde{\bm{z}},\bm{q}) \overline{\bm{p}_k\cdot \bm{u}^1_{\alpha,{\rm{tot}}}(\bm{x},\bm{z},\bm{q}_l) } \right)=\mathrm{Re}\left(\bm{p}_k\cdot \bm{u}^2_{\alpha,{\rm{tot}}}(\bm{x},\tilde{\bm{z}},\bm{q}) \overline{\bm{p}_k\cdot \bm{u}^2_{\alpha,{\rm{tot}}}(\bm{x},\bm{z},\bm{q}_l) } \right).
\end{equation}
According to \eqref{eq:uni1} and \eqref{eq:uni2}, we set
\begin{equation}\label{eq:thmdef}
\begin{aligned}
\bm{p}_k\cdot \bm{u}^j_{\alpha,{\rm{tot}}}(\bm{x},\tilde{\bm{z}},\bm{q}) &= r_{\alpha}^k(\bm{x},\tilde{\bm{z}},\bm{q}) {\rm{e}}^{{\rm{i}}\theta_{\alpha}^j(\bm{x},\tilde{\bm{z}},\bm{q}, \bm{p}_k) },\\
\bm{p}_k\cdot \bm{u}^j_{\alpha,{\rm{tot}}}(\bm{x},\bm{z},\bm{q}_l) &= s_{\alpha}^k(\bm{x},\bm{z},\bm{q}_l) {\rm{e}}^{{\rm{i}}\phi_{\alpha}^j(\bm{x},\bm{z},\bm{q}_l, \bm{p}_k) }, \quad j=1,2,
\end{aligned}
\end{equation}
where $r_{\alpha}^k(\bm{x},\tilde{\bm{z}},\bm{q})=|\bm{p}_k\cdot \bm{u}^j_{\alpha,{\rm{tot}}}(\bm{x},\tilde{\bm{z}},\bm{q}) |$, $s_{\alpha}^k(\bm{x},\bm{z},\bm{q}_l) = |\bm{p}_k\cdot \bm{u}^j_{\alpha,{\rm{tot}}}(\bm{x},\bm{z},\bm{q}_l) |$, and both $\theta_{\alpha}^j(\bm{x},\tilde{\bm{z}},\bm{q}, \bm{p}_k)$ and $\phi_{\alpha}^j(\bm{x},\bm{z},\bm{q}_l, \bm{p}_k)$ are real-valued functions.

To proceed to further proof, we prove the following lemma first.

\begin{lem}\label{lem:2d}
For all $\bm{z}\in\Sigma_1$, any fixed $l$, any fixed open set $\overline{\Gamma}_{h,1}\subset\Gamma_{h,1}$ and any $n\geq 2$ distinct $k_1, \cdots, k_n\in \{1,2,3,4,5\}$, at most one of $s_{\alpha}^{k_i}(\bm{x},\bm{z},\bm{q}_l)$ is identically zero for all $\bm{x}\in\overline{\Gamma}_{h,1}$.
\end{lem}
\begin{proof}\renewcommand{\qedsymbol}{}
 Otherwise, assume $s_{\alpha}^{k_1}(\bm{x}, \bm{z}_0, \bm{q}_l)$ and $s_{\alpha}^{k_2}(\bm{x}, \bm{z}_0, \bm{q}_l)$ are 0 for a fixed $\bm{z}_0\in \Sigma_1$ and all $\bm{x}\in \overline{\Gamma}_{h,1}$. Since $\bm{p}_1$ and $\bm{p}_2$ are not colinear as assumed, it follows that $\bm{u}^j_{\alpha,{\rm{tot}}}(\bm{x},\bm{z}_0,\bm{q}_l) =0, \forall \bm{x}\in\overline{\Gamma}_{h,1}$, $j=1,2$. From the analyticity, $\bm{u}^j_{\alpha,\rm{tot}}(\bm{x},\bm{z}_0,\bm{q}_l)=0, \forall \bm{x}\in\Gamma_{h,1}$. Due to the uniqueness of the direct scattering problem, we have $\bm{u}^j_{\alpha,{\rm{tot}}}(\bm{x},\bm{z}_0,\bm{q}_l) =0$ for $\bm{x}\in \{(x_1,x_2):0\leq x_1\leq 1, x_2>h \}$. Then the analyticity of $\bm{u}^j_{\alpha,{\rm{tot}}}(\bm{x},\bm{z}_0,\bm{q}_l) $ with respect to $\bm{x}$ implies $\bm{u}^j_{\alpha,{\rm{tot}}}(\bm{x},\bm{z}_0,\bm{q}_l) =0, \forall \bm{x}\in (\Omega_{1,1}^+\cap\Omega_{2,1}^+)\setminus \{\bm{z}_0 \}$. Note that $\bm{u}^j_{\alpha,{\rm{tot}}}(\bm{x},\bm{z}_0,\bm{q}_l) =\bm{u}^j_{\alpha,{\rm{sc}}}(\bm{x},\bm{z}_0,\bm{q}_l) +\bm{G}^{\alpha}(\bm{x},\bm{z}_0)\bm{q}_l $, from which we can conclude $|\bm{u}^j_{\alpha,{\rm{sc}}}(\bm{x},\bm{z}_0,\bm{q}_l)|=|\bm{G}^{\alpha}(\bm{x},\bm{z}_0)\bm{q}_l|\to\infty$ as $\bm{x}\to\bm{z}_0$. This leads to a contradiction. The same result still holds if we substitute one of $s_\alpha^{k_i}$ with $r_\alpha^{k_i}$.
\end{proof}

Without loss of generality, we assume that $s_\alpha^k(\bm{x},\bm{z},\bm{q}_1)$, $k=1,2,3,4$ are not identically zero for $\bm{z}\in\Sigma_1$ and $\bm{x}\in\Gamma_{h,1}$. Then there exists a point $(\bm{x}^*, \bm{z}^*)\in\Gamma_{h,1}\times \Sigma_1$ such that $s_\alpha^k(\bm{x}^*, \bm{z}^*, \bm{q}_1)\neq 0$. From the continuity, there exists a neighborhood $U_1\times V_1\subset \Gamma_{h,1}\times\Sigma_1$ of $(\bm{x}^*, \bm{z}^*)$ on which $s_\alpha^k(\bm{x},\bm{z},\bm{q}_1)\neq 0$. From the discussion above and again without loss of generality, we could further assume that $s_\alpha^k(\bm{x}, \bm{z}, \bm{q}_2)$, $k=1,2,3$ are not identically zero on $U_1\times V_1$, and there exists $U_2\times V_2\subset U_1\times V_1$ on which $s_\alpha^k(\bm{x},\bm{z}, \bm{q}_2)\neq 0$. Repeat the same steps again and we could finally find $U\times V\subset U_2\times V_2\subset\Gamma_{h,1}\times\Sigma_1$ such that $s_\alpha^k(\bm{x},\bm{z},\bm{q}_l)$ and $r_\alpha^k(\bm{x},\tilde{\bm{z}},\bm{q})$ are non-zero on $U\times V$ for $k=1,2$ and $l=1,2$.


For any $(\bm{x},\bm{z})\in U\times V$, \eqref{eq:realpart} implies
\[
{\rm{cos}}\left(\theta_{\alpha}^1(\bm{x},\tilde{\bm{z}},\bm{q}, \bm{p}_k)-\phi_{\alpha}^1(\bm{x},\bm{z},\bm{q}_l, \bm{p}_k) \right) = {\rm{cos}}\left(\theta_{\alpha}^2(\bm{x},\tilde{\bm{z}},\bm{q}, \bm{p}_k)-\phi_{\alpha}^2(\bm{x},\bm{z},\bm{q}_l, \bm{p}_k) \right),
\]
which means either
\begin{equation}\label{eq:angle1}
\theta_{\alpha}^1(\bm{x},\tilde{\bm{z}},\bm{q}, \bm{p}_k)-\theta_{\alpha}^2(\bm{x},\tilde{\bm{z}},\bm{q}, \bm{p}_k)= \phi_{\alpha}^1(\bm{x},\bm{z},\bm{q}_l, \bm{p}_k) -\phi_{\alpha}^2(\bm{x},\bm{z},\bm{q}_l, \bm{p}_k) +2n\pi
\end{equation}
or
\begin{equation}\label{eq:angle2}
\theta_{\alpha}^1(\bm{x},\tilde{\bm{z}},\bm{q}, \bm{p}_k)+\theta_{\alpha}^2(\bm{x},\tilde{\bm{z}},\bm{q}, \bm{p}_k)= \phi_{\alpha}^1(\bm{x},\bm{z},\bm{q}_l, \bm{p}_k) +\phi_{\alpha}^2(\bm{x},\bm{z},\bm{q}_l, \bm{p}_k) +2n\pi
\end{equation}
for $n\in\mathbb{Z}$, $k,l=1,2$.

We first consider the case \eqref{eq:angle1}. Define $\gamma_{\alpha}(\bm{x},\tilde{\bm{z}},\bm{q}, \bm{p}_k )=\theta_{\alpha}^1(\bm{x},\tilde{\bm{z}},\bm{q}, \bm{p}_k)-\theta_{\alpha}^2(\bm{x},\tilde{\bm{z}},\bm{q}, \bm{p}_k)-2n\pi$, we have
\begin{align}\label{eq:rela1}
\nonumber \bm{p}_k\cdot \bm{u}^1_{\alpha,{\rm{tot}}}(\bm{x},\bm{z},\bm{q}_l)= &s_{\alpha}^k(\bm{x},\bm{z},\bm{q}_l) {\rm{e}}^{{\rm{i}}\phi_{\alpha}^1(\bm{x},\bm{z},\bm{q}_l, \bm{p}_k) }\\
\nonumber =&s_{\alpha}^k(\bm{x},\bm{z},\bm{q}_l) {\rm{e}}^{{\rm{i}}\phi_{\alpha}^2(\bm{x},\bm{z},\bm{q}_l, \bm{p}_k) }  {\rm{e}}^{{\rm{i}} \gamma_{\alpha}(\bm{x}, \tilde{\bm{z}},\bm{q}, \bm{p}_k)}\\
=&\bm{p}_k\cdot \bm{u}^2_{\alpha,{\rm{tot}}}(\bm{x},\bm{z},\bm{q}_l){\rm{e}}^{{\rm{i}} \gamma_{\alpha}(\bm{x}, \tilde{\bm{z}},\bm{q}, \bm{p}_k)}, \quad \forall (\bm{x},\bm{z})\in U\times V.
\end{align}
According to the reciprocity relation in Lemma~\ref{lem:33}, \eqref{eq:rela1} shows that
\[
\bm{q}_l \cdot \bm{u}^1_{\hat{\alpha},{\rm{tot}}}(\bm{z},\bm{x},\bm{p}_k) = \bm{q}_l \cdot \bm{u}^2_{\hat{\alpha},{\rm{tot}}}(\bm{z},\bm{x},\bm{p}_k){\rm{e}}^{{\rm{i}} \gamma_{\alpha}(\bm{x}, \tilde{\bm{z}},\bm{q}, \bm{p}_k)}, \quad \forall (\bm{x},\bm{z})\in U\times V.
\]
In view of the linear independence of $\bm{q}_1$ and $\bm{q}_2$, we have
\[
\bm{u}^1_{\hat{\alpha},{\rm{tot}}}(\bm{z},\bm{x},\bm{p}_k) =\bm{u}^2_{\hat{\alpha},{\rm{tot}}}(\bm{z},\bm{x},\bm{p}_k){\rm{e}}^{{\rm{i}} \gamma_{\alpha}(\bm{x}, \tilde{\bm{z}},\bm{q}, \bm{p}_k)}, \quad \forall (\bm{x},\bm{z})\in U\times V,\ k=1,2.
\]
Let $w_{\hat{\alpha}}(\bm{z},\bm{x},\bm{p}_k)=\bm{u}^1_{\hat{\alpha},{\rm{tot}}}(\bm{z},\bm{x},\bm{p}_k) -\bm{u}^2_{\hat{\alpha},{\rm{tot}}}(\bm{z},\bm{x},\bm{p}_k){\rm{e}}^{{\rm{i}} \gamma_{\alpha}(\bm{x}, \tilde{\bm{z}},\bm{q}, \bm{p}_k)}$. Due to the analyticity of $w_{\hat{\alpha}}(\bm{z},\bm{x},\bm{p}_k)$ with respect to $\bm{z}$ and the fact that $V$ is an admissible curve, we obtain that for any $\bm{x}\in U$,
\[
w_{\hat{\alpha}}(\bm{z},\bm{x},\bm{p}_k)=0,\quad \forall\bm{z}\in \partial\Omega_1,
\]
where $\Omega_1$ is defined in Definition \ref{def:admissiblecurve}, hence
\begin{align}
\left\{
\begin{array}{ll}
\Delta^{\ast}w_{\hat{\alpha}}+\omega^2 w_{\hat{\alpha}}=0 &{\rm{in}} \ \Omega_1,\\
w_{\hat{\alpha}}=0 & {\rm{on}}\  \partial \Omega_1.
\end{array}
\right.
\end{align}
Since $\omega^2$ is not a Dirichlet eigenvalue of $-\Delta^{\ast}$ in $\Omega_1$, we obtain that for any $\bm{x}\in U$,
\[
w_{\hat{\alpha}}(\bm{z},\bm{x},\bm{p}_k)=0,\quad \forall\bm{z}\in \Omega_1.
\]
Then, the analyticity of $\bm{u}^j_{\hat{\alpha},{\rm{tot}}}(\bm{z},\bm{x},\bm{p}_k)$ with respect to $\bm{z}$ yields for $k=1,2$ that
\begin{equation}\label{eq:sc_inc}
\bm{u}^1_{\hat{\alpha},{\rm{sc}}}(\bm{z},\bm{x},\bm{p}_k) =\bm{u}^2_{\hat{\alpha},{\rm{sc}}}(\bm{z},\bm{x},\bm{p}_k) {\rm{e}}^{{\rm{i}} \gamma_{\alpha}(\bm{x}, \tilde{\bm{z}},\bm{q}, \bm{p}_k)} + \bm{G}_{\hat{\alpha}}(\bm{z},\bm{x})\bm{p}_k \left({\rm{e}}^{{\rm{i}} \gamma_{\alpha}(\bm{x}, \tilde{\bm{z}},\bm{q}, \bm{p}_k)}-1 \right)
\end{equation}
for any $\bm{z}\in (\Omega_{1,1}^+\cap \Omega_{2,1}^+)\setminus \{\bm{x} \}$. The boundedness of $\bm{u}^j_{\hat{\alpha},{\rm{sc}}}(\bm{z},\bm{x},\bm{p}_k)$ and the unboundedness of $\bm{G}_{\hat{\alpha}}(\bm{z},\bm{x})\bm{p}_k$ imply that \eqref{eq:sc_inc} holds if and only if
\[
{\rm{e}}^{{\rm{i}} \gamma_{\alpha}(\bm{x}, \tilde{\bm{z}},\bm{q}, \bm{p}_k)} =1, \quad \forall \bm{x}\in U.
\]
Therefore, one has
\[
\bm{u}^1_{\hat{\alpha},{\rm{tot}}}(\bm{z},\bm{x},\bm{p}_k) =\bm{u}^2_{\hat{\alpha},{\rm{tot}}}(\bm{z},\bm{x},\bm{p}_k),\  \forall (\bm{z},\bm{x})\in (\Omega_{1,1}^+\cap \Omega_{2,1}^+)\times U, \ \bm{z}\neq \bm{x}.
\]
By taking the dot product of the above equation with $\bm{q}_l\in \mathbb{S}$ and using the reciprocity relation in Lemma~\ref{lem:33}, we have
\[
\bm{p}_k\cdot\bm{u}^1_{\alpha,{\rm{tot}}}(\bm{x},\bm{z},\bm{q}_l) = \bm{p}_k\cdot\bm{u}^2_{\alpha,{\rm{tot}}}(\bm{x},\bm{z},\bm{q}_l), \quad \forall (\bm{x},\bm{z})\in U\times (\Omega_{1,1}^+\cap \Omega_{2,1}^+), \ \bm{x}\neq \bm{z},\ k,l=1,2.
\]
Since $\bm{p}_1$ and $\bm{p}_2$ are linearly independent, we get
\[
\bm{u}^1_{\alpha,{\rm{tot}}}(\bm{x},\bm{z},\bm{q}_l) =\bm{u}^2_{\alpha,{\rm{tot}}}(\bm{x},\bm{z},\bm{q}_l), \quad \forall (\bm{x},\bm{z})\in U\times (\Omega_{1,1}^+\cap \Omega_{2,1}^+), \ \bm{x}\neq \bm{z},\ l=1,2.
\]
By the analyticity of $\bm{u}^j_{\alpha,{\rm{tot}}}(\bm{x},\bm{z},\bm{q}_l)$ with respect to $\bm{x}$, we obtain that
\[
\bm{u}^1_{\alpha,{\rm{tot}}}(\bm{x},\bm{z},\bm{q}_l) =\bm{u}^2_{\alpha,{\rm{tot}}}(\bm{x},\bm{z},\bm{q}_l), \quad \forall (\bm{x},\bm{z})\in \Gamma_{h,1}\times (\Omega_{1,1}^+\cap \Omega_{2,1}^+),\ \bm{x}\neq \bm{z},\ l=1,2.
\]
Then, we get $f_1=f_2$ by applying the uniqueness result (cf. theorem 4 in \cite{Charalambopoulos01}).

Next, we shall show that the case \eqref{eq:angle2} does not hold. Suppose \eqref{eq:angle2} is true, then following the similar arguments as above, we observe that there exists a real-valued function $\eta_{\alpha}( \bm{x},\tilde{\bm{z}},\bm{q}, \bm{p}_k )$ such that
\[
\bm{u}^1_{\hat{\alpha},{\rm{tot}}}(\bm{z},\bm{x},\bm{p}_k) =\overline{\bm{u}^2_{\hat{\alpha},{\rm{tot}}}(\bm{z},\bm{x},\bm{p}_k)} {\rm{e}}^{{\rm{i}} \eta_{\alpha}(\bm{x}, \tilde{\bm{z}},\bm{q}, \bm{p}_k)}, \quad \forall (\bm{x},\bm{z})\in U\times (\Omega_{1,1}^+\cap \Omega_{2,1}^+),\ \bm{x}\neq \bm{z},\ k=1,2.
\]
Thus, the boundedness of $\bm{u}^j_{\hat{\alpha},{\rm{sc}}}(\bm{z},\bm{x},\bm{p}_k) $ shows that $\bm{G}_{\hat{\alpha}}(\bm{z},\bm{x})\bm{p}_k -  {\rm{e}}^{{\rm{i}} \eta_{\alpha}(\bm{x}, \tilde{\bm{z}},\bm{q}, \bm{p}_k)}\overline{\bm{G}_{\hat{\alpha}}(\bm{z},\bm{x})\bm{p}_k }$ is bounded for all $\bm{z}\in (\Omega_{1,1}^+\cap \Omega_{2,1}^+)\setminus \{\bm{x}\}$. Letting $\bm{z}\to \bm{x}$ along the straight line $z_1-x_1=z_2-x_2$, by the singularity of $\bm{G}_{\hat{\alpha}}(\bm{z},\bm{x})$, we have that $\mathrm{Im}(\bm{G}_{\hat{\alpha}}(\bm{z},\bm{x})\bm{p}_k)$ is bounded, and $\mathrm{Re}(\bm{G}_{\hat{\alpha}}(\bm{z},\bm{x})\bm{p}_k)$ is unbounded, which implies that $ {\rm{e}}^{{\rm{i}} \eta_{\alpha}(\bm{x}, \tilde{\bm{z}},\bm{q}, \bm{p}_k)}\equiv 1$. Consequently, we have
\[
\bm{u}^1_{\hat{\alpha},{\rm{tot}}}(\bm{z},\bm{x},\bm{p}_k) =\overline{\bm{u}^2_{\hat{\alpha},{\rm{tot}}}(\bm{z},\bm{x},\bm{p}_k)}
\]
for $\forall (\bm{x},\bm{z})\in U\times (\Omega_{1,1}^+\cap \Omega_{2,1}^+),\ \bm{x}\neq \bm{z},\ k=1,2$. Next, by taking the dot product with a fixed $\bm{q}$ and using Lemma~\ref{lem:33}, along with the linear independence of vectors $\bm{p}_1$ and $\bm{p}_2$, we have
\[
\bm{u}^1_{{\alpha},{\rm{tot}}}(\bm{x},\bm{z},\bm{q}) =\overline{\bm{u}^2_{{\alpha},{\rm{tot}}}(\bm{x},\bm{z},\bm{q})}, \quad \forall (\bm{x},\bm{z})\in U\times (\Omega_{1,1}^+\cap \Omega_{2,1}^+),\ \bm{x}\neq \bm{z}.
\]
Meanwhile, for any $\bm{x}\in U$, both $\bm{u}^1_{{\alpha},{\rm{tot}}}(\bm{x},\bm{z},\bm{q})$ and $\bm{u}^2_{{\alpha},{\rm{tot}}}(\bm{x},\bm{z},\bm{q})$ admit the radiating expansion \eqref{eq:sc_radiation}. By direct calculation, we observe that $\bm{u}^1_{{\alpha},{\rm{tot}}}(\bm{x},\bm{z},\bm{q})$ only contains upgoing waves and evanescent waves, while $\overline{\bm{u}^2_{{\alpha},{\rm{tot}}}(\bm{x},\bm{z},\bm{q})}$ only contains downgoing waves and evanescent waves. This implies that $\bm{u}^1_{{\alpha},{\rm{tot}}}(\bm{x},\bm{z},\bm{q})$ only contains evanescent waves, which contradicts Assumption~\ref{asm:1}.
\end{proof}

\subsection{Uniqueness of quasi-periodic/biperiodic elastic system in $\mathbb{R}^3$}
We will prove the uniqueness of the quasi-periodic/biperiodic elastic system in $\mathbb{R}^3$ in this section. Since the proof is very similar to Theorem \ref{thm:2d}, we will only demonstrate the different parts to keep the paper concise.
\begin{thm}
Let $f_2, \tilde{f}_2\in C^2(\mathbb{R})$ be 1-periodic functions. Suppose that the $\alpha$-quasi-periodic total fields satisfy
\begin{align*}
|\bm{p}_k\cdot \bm{u}^1_{\alpha,{\rm{tot}}}(\bm{x},\tilde{\bm{z}},\bm{q})| =& |\bm{p}_k\cdot \bm{u}^2_{\alpha,{\rm{tot}}}(\bm{x},\tilde{\bm{z}},\bm{q})|, \quad \forall \bm{x} \in \overline{\Gamma}_{h,2},\\
 |\bm{p}_k\cdot \bm{u}^1_{\alpha,{\rm{tot}}}(\bm{x},\bm{z},\bm{q}_l)| = & |\bm{p}_k\cdot \bm{u}^2_{\alpha,{\rm{tot}}}(\bm{x},\bm{z},\bm{q}_l)|, \quad \forall(\bm{x},\bm{z}) \in \overline{\Gamma}_{h,2}\times \Sigma_2,\\
 |\bm{p}_k\cdot \left(\bm{u}^1_{\alpha,{\rm{tot}}}(\bm{x},\tilde{\bm{z}},\bm{q})+ \bm{u}^1_{\alpha,{\rm{tot}}}(\bm{x},\bm{z},\bm{q}_l) \right)|  =& |\bm{p}_k\cdot \left(\bm{u}^2_{\alpha,{\rm{tot}}}(\bm{x},\tilde{\bm{z}},\bm{q})+ \bm{u}^2_{\alpha,{\rm{tot}}}(\bm{x},\bm{z},\bm{q}_l) \right)|,  \ \forall(\bm{x}, \bm{z}) \in \overline{\Gamma}_{h,2}\times\Sigma_2,
\end{align*}
for any fixed $\tilde{\bm{z}}\in\Gamma_{0,2}$, open subset $\overline{\Gamma}_{h,2}\subset\Gamma_{h,2}$, $\bm{q}, \bm{q}_l\in\mathbb{S}$, $l=1,2,3$, $k=1,2,\cdots,9$ and the wavenumbers $k_p\in [k_{\mathrm{min}}, k_{\mathrm{max}}]$ for some $0<k_{\mathrm{min}}<k_{\mathrm{max}}$. Here, $\Sigma_2$ is the admissible surface with respect to $\Omega_2$, and $\bm{q}_1,\bm{q}_2,\bm{q}_3$ and $\bm{p}_1,\bm{p}_2,\cdots,\bm{p}_9$ are any given polarization vectors that $\bm{q}_1, \bm{q}_2, \bm{q}_3$ and $\bm{p}_{i_1}, \bm{p}_{i_2}, \bm{p}_{i_3}$ are linearly independent for any distinct $i_1, i_2, i_3$ ranging from 1 to 9. Then $f_2$ and $\tilde{f}    _2$ coincide.
\end{thm}
\begin{proof}
We define $r_{\alpha}^k(\bm{x},\tilde{\bm{z}},\bm{q}), s_{\alpha}^k(\bm{x},\bm{z},\bm{q}_l), \theta_{\alpha}^j(\bm{x},\tilde{\bm{z}},\bm{q}, \bm{p}_k)$ and $\phi_{\alpha}^j(\bm{x},\bm{z},\bm{q}_l, \bm{p}_k)$ as in \eqref{eq:thmdef}. Lemma \ref{lem:2d} is replaced as follows.
\begin{lem}\label{lem:3d1}
For all $\bm{z}\in\Sigma$, any fixed $l$, any fixed open set $\overline{\Gamma}_{h,2}\subset\Gamma_{h,2}$ and any $n\geq 3$ distinct $k_1, \cdots, k_n\in \{1,2,\cdots,9\}$, at most two of $s_{\alpha}^{k_i}(\bm{x},\bm{z},\bm{q}_l)$ is identically zero for all $\bm{x}\in\overline{\Gamma}_{h,2}$.
\end{lem}
\begin{proof}\renewcommand{\qedsymbol}{}
 We prove by contradiction. Assume $s_{\alpha}^{k_i}(\bm{x}, \bm{z}_0, \bm{q}_l), i=1,2,3$ are 0 for a fixed $\bm{z}_0\in \Sigma$ and all $\bm{x}\in \overline{\Gamma}_{h,2}$. Since $\bm{p}_i, i=1,2,3$ are linearly independent, we have $\bm{u}^j_{\alpha,{\rm{tot}}}(\bm{x},\bm{z}_0,\bm{q}_l) =0, \forall \bm{x}\in\overline{\Gamma}_{h,2}$, $j=1,2$. The rest part is the same as proof of Lemma \ref{lem:2d}.
\end{proof}
Now without loss of generality, we assume that $s_\alpha^k(\bm{x},\bm{z},\bm{q}_1)$, $k=1,\cdots,7$ are not identically zero for $\bm{z}\in\Sigma_2$ and $\bm{x}\in\Gamma_{h,2}$. Then there exists a point $(\bm{x}^*, \bm{z}^*)\in\Gamma_{h,2}\times \Sigma_2$ such that $s_\alpha^k(\bm{x}^*, \bm{z}^*, \bm{q}_1)\neq 0$. From the continuity, there exists a neighborhood $U_1\times V_1\subset \Gamma_{h,2}\times\Sigma_2$ of $(\bm{x}^*, \bm{z}^*)$ on which $s_\alpha^k(\bm{x},\bm{z},\bm{q}_1)\neq 0$. From the discussion above and again without loss of generality, we could further assume that $s_\alpha^k(\bm{x}, \bm{z}, \bm{q}_2)$, $k=1,\cdots,5$ are not identically zero on $U_1\times V_1$, and there exists $U_2\times V_2\subset U_1\times V_1$ on which $s_\alpha^k(\bm{x},\bm{z}, \bm{q}_2)\neq 0$. Repeat the same steps again and we could finally find $U\times V\subset U_2\times V_2\subset\Gamma_{h,2}\times\Sigma_2$ such that $s_\alpha^k(\bm{x},\bm{z},\bm{q}_l)$ and $r_\alpha^k(\bm{x},\tilde{\bm{z}},\bm{q})$ is non-zero on $U\times V$ for $k=1,2,3$ and $l=1,2,3$. The remaining part is the same as proof of Theorem \ref{thm:2d}, except we apply Theorem \ref{thm:unique1} in the appendix for the last step of uniqueness proof.
\end{proof}

\begin{thm}
Let $f_3, \tilde{f}_3\in C^2(\mathbb{R}^2)$ be biperiodic functions with periodicity $(1,1)$. Suppose that the $\bm{\alpha}$-quasi-biperiodic total fields satisfy
\begin{align*}
|\bm{p}_k\cdot \bm{u}^1_{\bm{\alpha},{\rm{tot}}}(\bm{x},\tilde{\bm{z}},\bm{q})| =& |\bm{p}_k\cdot \bm{u}^2_{\bm{\alpha},{\rm{tot}}}(\bm{x},\tilde{\bm{z}},\bm{q})|, \quad \forall \bm{x} \in \overline{\Gamma}_{h,3},\\
 |\bm{p}_k\cdot \bm{u}^1_{\bm{\alpha},{\rm{tot}}}(\bm{x},\bm{z},\bm{q}_l)| = & |\bm{p}_k\cdot \bm{u}^2_{\bm{\alpha},{\rm{tot}}}(\bm{x},\bm{z},\bm{q}_l)|, \quad \forall(\bm{x},\bm{z}) \in \overline{\Gamma}_{h,3}\times \Sigma_3,\\
 |\bm{p}_k\cdot \left(\bm{u}^1_{\bm{\alpha},{\rm{tot}}}(\bm{x},\tilde{\bm{z}},\bm{q})+ \bm{u}^1_{\bm{\alpha},{\rm{tot}}}(\bm{x},\bm{z},\bm{q}_l) \right)|  =& |\bm{p}_k\cdot \left(\bm{u}^2_{\bm{\alpha},{\rm{tot}}}(\bm{x},\tilde{\bm{z}},\bm{q})+ \bm{u}^2_{\bm{\alpha},{\rm{tot}}}(\bm{x},\bm{z},\bm{q}_l) \right)|,  \ \forall(\bm{x}, \bm{z}) \in \overline{\Gamma}_{h,3}\times\Sigma_3,
\end{align*}
for any fixed $\tilde{\bm{z}}\in\Gamma_{0,3}$, open subset $\overline{\Gamma}_{h,3}\subset\Gamma_{h,3}$, $\bm{q}, \bm{q}_l\in\mathbb{S}$, $l=1,2,3$, $k=1,2,\cdots,9$ and the wavenumbers $k_p\in [k_{\mathrm{min}}, k_{\mathrm{max}}]$ for some $0<k_{\mathrm{min}}<k_{\mathrm{max}}$. Here, $\Sigma_3$ is the admissible surface with respect to $\Omega_3$, and $\bm{q}_1,\bm{q}_2,\bm{q}_3$ and $\bm{p}_1,\bm{p}_2,\cdots,\bm{p}_9$ are any given polarization vectors that $\bm{q}_1, \bm{q}_2, \bm{q}_3$ and $\bm{p}_{i_1}, \bm{p}_{i_2}, \bm{p}_{i_3}$ are linearly independent for any distinct $i_1, i_2, i_3$ ranging from 1 to 9. Then $f_3$ and $\tilde{f}_3$ coincide.
\end{thm}

\section{Acknowledgements}
Y. H. was supported by NSFC grant (12401564) and Shenzhen Science and Technology Program (Grant No. RCBS20231211090750090). W. W. was supported by NSFC grant (12301539) and youth growth grant of Department of Science and Technology of Jilin Province (20240602113RC). 

\bibliographystyle{plain}
\bibliography{unique}

\begin{appendix}
\section{Uniqueness theorem in three-dimensional space}

In this section we would like to prove the uniqueness theorem in quasi-periodic three-dimensional system, which will be needed in following sections. The proof follows the proof of uniqueness theorem in two-dimensional space in \cite{Charalambopoulos01}. Since the steps only differs from \cite{Charalambopoulos01} in some details, here we will only point out the different places instead of giving a full proof to the lemmas and theorem. The proof of uniqueness theorem in quasi-biperiodic three-dimensional system is exactly the same as quasi-periodic three-dimensional one, except with a different periodic unit, so we will not repeat the proof here.

Firstly we need a few necessary lemmas. We consider the layer
$$
    X:=\{\bm{x}\in\mathbb{R}^3: f(x_1) < x_3 < g(x_1), x_2\in\mathbb{R}\}
$$
where $f,g\in C^2(\mathbb{R})$ denote two 1-periodic functions with $f(t)<g(t)$ for all $t\in\mathbb{R}$.

\begin{lem}\label{lem:lem1}
Let $k_p:=\omega/c_p\in\mathbb{C}$, where $\omega$ denotes frequency and $c_p:=\sqrt{(\lambda+2\mu)/\rho}$ denotes the speed of p-wave. Assume that $\mathrm{Im} k_p\neq 0$ and $\bm{u}\in C^2(X)\cap C(\overline{X})$ is a quasi-periodic function with quasi-periodicity $\alpha$ satisfying \eqref{eq:Navier} in $X$ and vanishing on $\partial X$. Then $\bm{u}$ vanishes in $X$. The same result holds if $k_p=0$.
\end{lem}

The proof of Lemma \ref{lem:lem1} is the same as Lemma 1 of \cite{Charalambopoulos01}, except we use
$$
    \mathcal{E}(\bm{v}, \bm{w}) = \lambda(\nabla\cdot\bm{v})(\nabla\cdot\bm{w}) + \mu\sum\limits_{i,j=1}^3\frac{\partial v_i}{\partial x_j}\left(\frac{\partial w_i}{\partial x_j} + \frac{\partial w_j}{\partial x_i}\right)
$$
as the quadratic form.

\begin{defn}
Let $m\in (-1,1)$ be fixed. The wavenumber $k_p \in\mathbb{R}$, $k_p \neq 0$ is called a Dirichlet eigenvalue of the layer $X$ if there exists a quasi-periodic nontrivial solution $\bm{u} \in [C^2(X) \cap C(\overline{X})]$ with momentum $k_pm$ of the equation of linearized elasticity with boundary condition $\bm{u} = 0$ on $\partial X$. The function $\bm{u}$ is called the corresponding eigenfunction.
\end{defn}

Let $\Omega\subset\mathbb{R}^3$ satisfying the cone condition, $L^2(\Omega,\Delta^{\ast}):=\{\bm{u}\in L^2(\Omega) : \Delta^{\ast}\bm{u}\in L^2(\Omega)\}$. We call $\bm{u}$ a weak solution of Dirichlet problem of Lam\'{e} equation in $\Omega$, if it belongs to $H^1(\Omega)\cap L^2(\Omega,\Delta^{\ast})$ and satisfies the weak formulation
\begin{equation}\label{eq:lameweak}
    \int_\Omega \lambda(\nabla\cdot\bm{u})(\nabla\cdot\overline{\bm{v}}) + 2\mu\epsilon_{ij}(\bm{u})\epsilon_{ij}(\overline{\bm{v}}) - \rho\omega^2\bm{u}\cdot\overline{\bm{v}} \mathrm{d}\bm{x} = 0
\end{equation}
for any $\bm{v}\in H^1(\Omega)$. In the quasi-periodic case, however, we need to slightly adjust the previous weak formulation. For quasi-periodic functions $\bm{u}$ and $\bm{v}$, denote $\bm{u}(\bm{x}) = \e^{\ii kmx_1}\bm{\phi}(\bm{x})$, $\bm{v}(\bm{x}) = \e^{\ii kmx_1}\bm{\psi}(\bm{x})$. Then $\bm{\phi}$ and $\bm{\psi}$ are both 1-periodic functions with respect to $x_1$. Define
\begin{align*}
    &M:=\{ \bm{x}\in\mathbb{R}^3: f(x_1)<x_3<g(x_1), 0<x_1<1, x_2\in\mathbb{R}\}, \\
    &N:=\{\bm{x}\in\mathbb{R}^3: x_3>f(x_1), x_1, x_2\in\mathbb{R}\}.
\end{align*}

Define Hilbert space $H$ as the completion of $\{\phi\in C_p^1(\overline{M})\cap H^1(M):\phi=0~ \mathrm{on}~ \partial N\cap M\}$ with respect to $H^1$-norm. The space $C_p^1(\overline{M})$ contains all the differentiable functions which are $1-$periodic with respect to $x_1$. If we work on periodic cell $M$, weak formulation \eqref{eq:lameweak} will be reformulated as
\begin{equation}\label{eq:lameweakperiodic}
\begin{aligned}
    &\int_M 2\mu\nabla\bm{\phi}:\nabla\overline{\bm{\psi}}^\top + \lambda(\nabla\cdot\bm{\phi})(\nabla\cdot\overline{\bm{\psi}}) + \mu(\nabla\times\bm{\phi})\cdot(\nabla\times\overline{\bm{\psi}})+ \ii k_p m\left[\lambda\bm{e}_1\cdot(\bm{\phi}(\nabla\cdot\overline{\bm{\psi}})- \overline{\bm{\psi}}(\nabla\cdot\bm{\phi}))\right. \\
    & \left.  + (\bm{e}_3\otimes\bm{e}_2 - \bm{e}_2\otimes\bm{e}_3):\mu((\nabla\times\bm{\phi})\otimes\overline{\bm{\psi}} - (\nabla\times\overline{\bm{\psi}})\otimes\bm{\phi}) +\bm{e}_1\cdot 2\mu(\nabla\overline{\bm{\psi}}\cdot\bm{\phi} - \nabla\bm{\phi}\cdot\overline{\bm{\psi}})\right]\\
    & + k_p^2\left[m^2((\lambda+2\mu)\bm{e}_1\otimes\bm{e}_1 + \mu(\bm{e}_2\otimes\bm{e}_2 + \bm{e}_3\otimes\bm{e}_3)):\bm{\phi}\otimes\overline{\bm{\psi}} - (\lambda+2\mu)\bm{\phi}\cdot\overline{\bm{\psi}}\right] \mathrm{d}\bm{x} = 0.
\end{aligned}
\end{equation}
Here the double contraction of dyads is defined as $\bm{a}\otimes\bm{b}:\bm{c}\otimes\bm{d} = (\bm{a}\cdot\bm{c})(\bm{b}\cdot\bm{d})$.

We equip $H$ with the inner product defined by
$$
    \langle \bm{\phi}, \bm{\psi}\rangle_H = \int_M 2\mu\nabla\bm{\phi}:\nabla\overline{\bm{\psi}}^\top + \lambda(\nabla\cdot\bm{\phi})(\nabla\cdot\overline{\bm{\psi}}) + \mu(\nabla\times\bm{\phi})\cdot(\nabla\times\overline{\bm{\psi}}) \mathrm{d}\bm{x}
$$
for $\bm{\phi}, \bm{\psi}\in H$.

Now we could prove the following lemma.

\begin{lem}
    For any $k_p^*>0$ and fixed $m\in (-1,1)$, the Dirichlet eigenvalues $|k_p|\leq k_p^*$ of the periodic layer $X$ form a finite set.
\end{lem}
The proof of this lemma is the same as Theorem 1 in \cite{Charalambopoulos01}, except we replace the sesquilinear forms $a_1(\bm{\phi}, \bm{\psi})$ and $a_2(\bm{\phi}, \bm{\psi})$ with
\begin{align*}
    a_1(\bm{\phi}, \bm{\psi}) &= \ii m\int_M \lambda\bm{e}_1\cdot(\bm{\phi}(\nabla\cdot\overline{\bm{\psi}})- \overline{\bm{\psi}}(\nabla\cdot\bm{\phi})) \\ &+ (\bm{e}_3\otimes\bm{e}_2 - \bm{e}_2\otimes\bm{e}_3):\mu((\nabla\times\bm{\phi})\otimes\overline{\bm{\psi}} - (\nabla\times\overline{\bm{\psi}})\otimes\bm{\phi}) \\ &+\bm{e}_1\cdot 2\mu(\nabla\overline{\bm{\psi}}\cdot\bm{\phi} - \nabla\bm{\phi}\cdot\overline{\bm{\psi}})\mathrm{d}\bm{x}, \\
    a_2(\bm{\phi}, \bm{\psi}) &= \int_M \left[m^2((\lambda+2\mu)\bm{e}_1\otimes\bm{e}_1 + \mu(\bm{e}_2\otimes\bm{e}_2 + \bm{e}_3\otimes\bm{e}_3)):\bm{\phi}\otimes\overline{\bm{\psi}} - (\lambda+2\mu)\bm{\phi}\cdot\overline{\bm{\psi}}\right] \mathrm{d}\bm{x}
\end{align*}

Finally we could prove the uniqueness theorem as follows.

\begin{thm}\label{thm:unique1}
Let $f,g\in C^2(\mathbb{R})$ be $1$-periodic functions and $a>\max\{f(t), g(t): t\in [0,1]\}$. Define surface $\Gamma_f:= \{(x_1,x_2,x_3):x_3=f(x_1), x_2\in\mathbb{R}\}$ and $\Gamma_g$ is defined similarly. We require that the incident wave is not perpendicular to $(x_1,x_2)$-plane.  Assuming that the wave scattered by $\Gamma_f$ and $\Gamma_g$ are equivalent, i.e.
$$
    \bm{u}_{sc,f}(\cdot,a) = \bm{u}_{sc,g}(\cdot,a)
$$
for all incident plane waves $\widetilde{\bm{G}}^\alpha(\bm{x},\bm{z})\bm{q}$ with frequency $\omega\in[\omega_{\mathrm{min}}, \omega_{\mathrm{\max}}]$, where $0<\omega_{\mathrm{min}}<\omega_{\mathrm{max}}$. Then $\Gamma_f$ and $\Gamma_g$ coincide.
\end{thm}

The uniqueness theorem for quasi-biperiodic system, which could be proved in the exact same way, is stated as follows.
\begin{thm}\label{thm:unique2}
Let $\tilde{f},\tilde{g}\in C^2(\mathbb{R}^2)$ be biperiodic functions with periodicity $(1,1)$, and $b>\max\{f(t), g(t): t\in [0,1]^2\}$. Define surface $\Gamma_{\tilde{f}}:= \{(x_1,x_2,x_3):x_3=\tilde{f}(x_1,x_2)\}$ and $\Gamma_{\tilde{g}}$ is defined similarly. We require that the incident wave is not perpendicular to $(x_1,x_2)$-plane. Assuming that the wave scattered by $\Gamma_{\tilde{f}}$ and $\Gamma_{\tilde{g}}$ are equivalent, i.e.
$$
    \bm{u}_{sc,\tilde{f}}(\cdot,b) = \bm{u}_{sc,\tilde{g}}(\cdot,b)
$$
for all incident plane waves $\mathring{\bm{G}}^\alpha(\bm{x},\bm{z})\bm{q}$ with frequency $\omega\in[\omega_{\mathrm{min}}, \omega_{\mathrm{\max}}]$, where $0<\omega_{\mathrm{min}}<\omega_{\mathrm{max}}$. Then $\Gamma_{\tilde{f}}$ and $\Gamma_{\tilde{g}}$ coincide.
\end{thm}
\end{appendix}

\end{document}